\documentstyle[12pt,psfig]{article}
\psfigurepath{ps}
\newfont{\ghot}{eufm10 scaled \magstep 2} % Large @ 12pt == 1440
\newfont{\sets}{msbm10 scaled \magstep 2} % Large @ 12pt == 1440

\def\be{\begin{equation}}
\def\ee{\end{equation}}
\def\bl{\rule[-1mm]{2.4mm}{2.4mm}}
\newtheorem{thrm}{T~h~e~o~r~e~m}
\newtheorem{lmm}{L~e~m~m~a}

\begin{document}
\title{Pictorial representation for antisymmetric eigenfunctions of  $PS-3$ integral equations}
\author{\copyright 2007 ~~~~A.B.~Bogatyrev
\thanks{Supported by RFBR grant 05-01-01027 and grant MD4798.2007.01}}
\date{}
\maketitle

\begin{center}
\parbox{14cm}
{\small

  Eigenvalue problem for Poincare-Steklov-3 integral equation is reduced to the
  solution of three transcendential equations for three unknown numbers, moduli
  of pants. The complete list of antisymmetric eigenfunctions of integral
  equation in terms of Kleinian membranes is given.
  }
\end{center}

\section{Introduction}

  Traditionally, integral equations are the
  subject of functional analysis and operator theory. In the
  contrast we show that methods of complex
  geometry and combinatorics are efficient for the study of
  the following singular integral Poincare-Steklov (briefly, PS) equations

\begin{equation}
\label{PSE}
\lambda~
V.p.\int\limits_I
\frac{u(t)}{t-x}\,dt -
V.p.\int\limits_I
\frac{u(t)\,dR(t)}{R(t)-R(x)}
=const,\quad x\in I:=(-1,~1),
\end{equation}
where $\lambda$ is the spectral parameter; $u(t)$ is the unknown function;
$const$ is independent of $x$. The functional parameter $R(t)$ of the equation
is a given smooth {\it nondegenerate} change of variable on the interval $I$:
\begin{equation}
0<\left|
\frac d{dt}
R(t)
\right|<\infty,
\qquad {\rm when}~ t\in [-1,1].
\label{ND}
\end{equation}
Under the assumption that $R(t)=:R_3(t)$ is a rational degree three function
with separate real critical values different from the endpoints of the interval $I$,
we give the constructive representation for the eigenvalues $\lambda$ and eigenfunctions
$u(x)$ of equation (\ref{PSE}). First we say a few words about the origin of PS 
integral equations and the related background.

\paragraph{Spectral Boundary Value Problem.}
Let a domain in the plane be subdivided into two simply connected
domains $\Omega_1$ and $\Omega_2$ by a smooth simple arc $\Gamma$.
We are looking for the values of the spectral parameter $\lambda$
when the following problem has nonzero solution:

{\it Find a harmonic function $U_s$ in the domain $\Omega_s$, $s=1,2$,
vanishing on the outer portion of the boundary:
$\partial\Omega_s\setminus\Gamma$. On the interface
$\Gamma$ the functions $U_1$ and $U_2$
coincide while their normal derivatives differ by the factor of
$-\lambda$:}
\be
\label{InterfaceJump}
-\lambda\frac{\partial U_1}{\partial n}=
\frac{\partial U_2}{\partial n}.
\ee

\paragraph{Applications.}

  Boundary value problems for the Laplace equation with spectral
  parameter in the boundary conditions were first considered by
  H.Poincare (1896) and V.A.Steklov (1901). The problems of this 
  kind  arise e.g. in the diffraction, (thermo-) conductivity of 
  composite materials and analysis of 2D models of  oil
  extraction.

  This particular
  problem (\ref{InterfaceJump}) arises in justification and
  optimization of {\it domain decomposition method} for the
  solution of boundary values problems for elliptic PDE. The
  eigenvalues $\lambda$ of the spectral problem and the traces of
  eigenfunctions $U_1=U_2$ on the interface $\Gamma$ are
  respectively the critical values and critical points of the
  following functional, the ratio of two Dirichlet integrals

\be
F(U)=\frac{\int_{\Omega_2}|\nabla U_2|^2 ~d\Omega_2}
{\int_{\Omega_1}|\nabla U_1|^2 ~d\Omega_1},
 \qquad U\in H_{oo}^{1/2}(\Gamma),
\label{Rayleigh}
\ee
  where $U_s$ is the harmonic continuation of the function $U$
  from interface $\Gamma$ to the domain $\Omega_s$, $s=1,2$,
  vanishing at the outer boundary of the domain.

\paragraph{Integral Equation.}
The reduction of the stated above boundary value problem to the interface
brings to the equation (\ref{PSE}). Let $V_s$ be the harmonic function conjugate
to $U_s$, $s=1,2$. From Cauchy-Riemann equations and (\ref{InterfaceJump}) it follows
that the tangent to the interface derivatives of $V_1$ and $V_2$ differ by the same factor
$-\lambda$. Integrating along $\Gamma$ we get

\begin{equation}
\lambda V_1(y)+V_2(y)=const, \quad y\in\Gamma. \label{HarmConj}
\end{equation}

For the half-plane the boundary values of conjugate harmonic
functions are related via Hilbert transform. To take advantage of
this transformation we consider a conformal mapping $\omega_s(y)$
from $\Omega_s$ to the open upper halfplane $\mbox{\sets H}$ with
normalization $\omega_s(\Gamma)=I$, $s=1,2$. Now equation
(\ref{HarmConj}) may be rewritten as

$$
-\frac{\lambda}{\pi}~
V.p.\int\limits_I
\frac{U_1(\omega_1^{-1}(t))}{t-\omega_1(y)}\,dt -
\frac1{\pi}~
V.p.\int\limits_I
\frac{U_2(\omega_2^{-1}(t'))}{t'-\omega_2(y)}\,dt'
=const,
\quad y\in\Gamma.
$$
Introducing new notation
$x:=\omega_1(y)\in I$;
$R:=\omega_2\circ\omega_1^{-1}:~~I\to\Gamma\to I$;
$u(t):=U_1(\omega_1^{-1}(t))$ and the change of variable
$t'=R(t)$ in the second integral, we arrive at the Poincare-Steklov equation
(\ref{PSE}). Note that here $R(t)$ is the decreasing function on $I$.

  Operator analysis of equivalent spectral problems, boundary
  value problem (\ref{InterfaceJump}) and Poincare-Steklov
  equation (\ref{PSE}), may be found e.g.in \cite{Bog0}. Here
  we only mention that the spectrum is discreet if (\ref{ND}) holds, the eigen values
  are positive and converge to $\lambda=1$.

\paragraph{Philosophy of the Research.}
  The aim of our study is to give explicit expressions for the eigen pairs
  $(\lambda,u)$ of the PS integral equation. For the rational degree two
  functions $R(x)=R_2(x)$ the eigen pairs  were expressed in terms of elliptic
  functions \cite{Bog3}. Next natural step is to consider degree three
  rational functions.

  Here the notion of {\it explicit solution} should be specified. Usually this
  term means an {\it elementary function} of parameters or a {\it quadrature} of
  it or application of other {\it permissible operations} (e.g. the concept of
  Umemura classical functions). The history of mathematics however knows many
  disappointing results when the solution of the prescribed form does not exist.
  The nature always forces us to introduce new types of transcendent objects to
  enlarge the scope of search. Cf.: "Mais cette \'etude intime de la nature des
  fonctions integrales ne peut se faire que par l'introduction de transcendantes
  nouvelles" \cite{Poi}.

  From the
  philosophical point of view our goal is to study the nature of
  the solutions of integral equation (\ref{PSE}) and the means for their 
  constructive representation.

  %Take for instance the algebraic equations. The
  %ancients were  able to solve quadratic equations. But after the
  %nvention  of the formulas for cubic  and quartic
  %equations  in the 16th century no progress was made
  %untill the 19th century when it became clear
  %that no formula including arithmetic operations and
  %radicals only can solve equations of the higher order. After
  %that Ch.Hermite and L.Kronecker suggested formulae involving
  %elliptic modular function to find the roots of degree five
  %equations. The ideas of C.Jordan elaborated by H.Umemura
  %resulted in a formula (involving hyperelliptic integrals and
  %theta constants) for the roots of arbitrary degree polynomial.
  %In modern mathematical physics it is very often that the
  %problems are "explicitly" solved in terms of suitable
  %transcendential functions, say solutions of Painleve equations.

  \paragraph{Brief Description of the Result.}
   Given rational degree three function $R_3(t)$, we explicitly associate
  it to a {\it pair of pants} in Sect. \ref{PantsOfR3}. On the
  other hand, given spectral parameter $\lambda$ and two
  auxiliary real parameters, we explicitly construct in Sect.
  \ref{Result} another pair of pants which additionally depend on
  one or two integers. When the above two pants are
  conformally equivalent, $\lambda$ is the eigenvalue of the
  PS integral equation with parameter $R_3(x)$. Essentially,
  this means that to find the spectrum of the given integral
  equation (\ref{PSE}) one has to solve three transcendential
  equations involving three {\it moduli of pants}.

  Whether this representation of the solutions may be considered as constructive 
  or not is a matter of discussion. On the one hand, today it is possible to 
  numerically evaluate the conformal structure of surfaces (e.g. via circle 
  packing). On the other, this representation allows us to obtain valuable 
  features of the solution: to find the number of zeroes of eigenfunction 
  $u(t)$, to localize the spectrum and to show the discrete mechanism of 
  generating the eigenvalues.

\section{Space of PS-3 Equations}

In what follows we consider integral equations (\ref{PSE}) with rational degree tree real
functional parameter $R(x)=R_3(x)$ and call such equations PS-3. We restrict ourselves to the case when
$R_3(x)$ has {\it four distinct real critical values different from} $\pm1$.
  The details of our further constructions depend on the
  topological properties of functional parameter of
  the integral equation. One may encounter one of five described in
  section \ref{Classification} typical situations ${\cal A}$, ${\cal B}1$, ${\cal B}21$,
  ${\cal B}22$, ${\cal B}23$ corresponding to the components in the space of
  admissible functions $R_3(x)$.

\subsection{Topology of the Branched Covering}
\label{TopBrCov}
  Degree three rational function $R_3(x)$ defines the three-
  sheeted branched covering of a Riemann sphere by another Riemann
  sphere. The Riemann-Hurwitz formula suggests that $R_3(x)$
  typically has four separate branch points $a_s$, $s=1,4$. This
  means that every value $a_s$ is covered by a critical point
  $b_s$, and an ordinary  point $c_s$. We have assumed that
  all four branch points $a_s$ are distinct, real and differ
  from $\pm 1$. Other possible configurations are discussed in
  \cite{Bog1}.

  Every point $y\neq a_s$ of the extended real axis
  $\hat{\mbox{\sets R}}:=\mbox{\sets R}\cup\{\infty\}$ belongs to
  exactly  one of two types. For the type (3:0) the pre-image
  $R_3^{-1}(y)$ consists of three distinct real points. For the
  type (1:2) the pre-image consists of a real and two complex
  conjugate points. The type of the point is locally constant on
  the extended real axis and changes when we step over the branch
  point. Let the branch points $a_s$ be enumerated in the
  natural cyclic order of $\hat{\mbox{\sets R}}$ so that the
  intervals $(a_1,a_2)$ and $(a_3,a_4)$ are filled with the
  points of the type (1:2). We specify the way to
  exclude the relabeling $a_1\leftrightarrow a_3$,
  $a_2\leftrightarrow a_4$ of branch points in Sect. \ref{Classification}.

    \label{LabelBranchPoint}

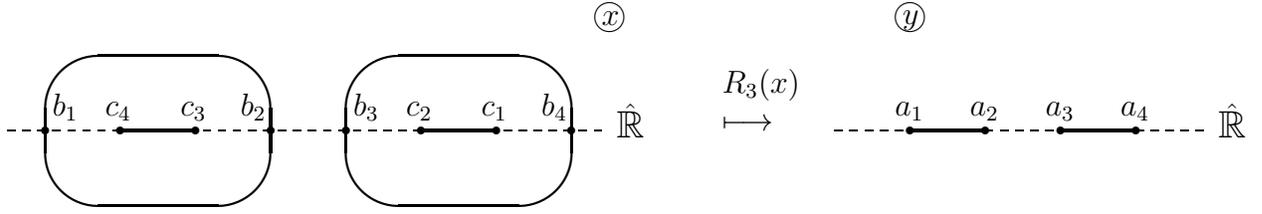
\begin{figure}
\begin{picture}(170,20)
%%%%%%%%%%%%%%%%%%%%%%%%%%%%%%%%%%%%%%%
%             lower sphere            %
%%%%%%%%%%%%%%%%%%%%%%%%%%%%%%%%%%%%%%%
\put(110,0){
\begin{picture}(70,20)

\thicklines
%cut D3
\put(10,10){\circle*{1}}
\put(20,10){\circle*{1}}
\put(10,10){\line(1,0){10}}
\put(8,12){$a_1$}
\put(18,12){$a_2$}
%cut D2
\put(30,10){\circle*{1}}
\put(40,10){\circle*{1}}
\put(30,10){\line(1,0){10}}
\put(28,12){$a_3$}
\put(38,12){$a_4$}

\thinlines
\multiput(0,10)(2,0){25}{\line(1,0){1}}
\put(51,9){$\hat{\mbox{\sets R}}$}
\put(10,25){\circle{4}}
\put(9,24.2){$y$}

\end{picture}
}
%%%%%%%%%%%%%%%%%%%%%%%%%%%%%%%%%%%%%%%
%             upper sphere            %
%%%%%%%%%%%%%%%%%%%%%%%%%%%%%%%%%%%%%%%
\put(5,0){
\begin{picture}(70,20)
%cuts
\thicklines
%cut R^-1(D2)
\put(10,10){\circle*{1}}
\put(20,10){\circle*{1}}
\put(10,10){\line(1,0){10}}
\put(8,12){$c_4$}
\put(18,12){$c_3$}
%cut R^-1(D1)
\put(50,10){\circle*{1}}
\put(60,10){\circle*{1}}
\put(50,10){\line(1,0){10}}
\put(48,12){$c_2$}
\put(58,12){$c_1$}
%ovals R^-1(D1,D2)
\put(15,10){\oval(30,20)}
\put(0,10){\circle*{1}}
\put(30,10){\circle*{1}}
\put(1,12){$b_1$}
\put(26,12){$b_2$}
\put(55,10){\oval(30,20)}
\put(40,10){\circle*{1}}
\put(70,10){\circle*{1}}
\put(41,12){$b_3$}
\put(66,12){$b_4$}
\thinlines
\multiput(-5,10)(2,0){40}{\line(1,0){1}}
\put(76,9){$\hat{\mbox{\sets R}}$}
\put(90,10){$\longmapsto$}
\put(90,15){$R_3(x)$}
\put(75,25){\circle{4}}
\put(74,24){$x$}

\end{picture}}
\end{picture}
\caption[]
{\small The topology of the covering $R_3$ with real branch points}
\label{CoverTopology}
\end{figure}

  The total pre-image $R_3^{-1}(\hat{\mbox{\sets R}})$ consists of
  the extended real axis  and two pairs of complex conjugate arcs
  intersecting $\hat{\mbox{\sets R}}$ at points $b_1$, $b_2$,
  $b_3$, $b_4$ as shown at the left side of Fig.
  \ref{CoverTopology}. The compliment of this pre-image on the
  Riemann sphere has six components, each of them is mapped 1-1
  onto upper or lower half plane. Note that the points $b_1$,
  $c_4$, $c_3$, $b_2$... on the left picture of Fig. \ref{CoverTopology}
  may follow in inverse order as well.

\subsection{Classification of Parameters $R_3$}
\label{Classification}

  The functional parameter $R_3(x)$ is a
  nondegenerate change of variable on the segment $[-1,1]$. This
  in particular means that no critical point $b_s$ belongs to this
  segment. So exactly one of two cases is realized:\footnote{
  Two points on a circle (extended real axis) define two segments.
  It should be clear which segment we mean: e.g. $b_1,b_2\not\in$ $[b_3,b_4]$;~~
  $b_1,b_4\not\in$ $[b_2,b_3]$, etc.}

\be
\begin{array}{cl}
Case~~{\cal A}: & [-1,1]\subset[b_2,b_3],\\
Case~~{\cal B}: & [-1,1]\subset[b_3,b_4].
\end{array}
\ee
  Remaining possibilities (like $[-1,1]\subset[b_1,b_2]$) are
  reduced either to ${\cal A}$ or ${\cal B}$ by the clever choice of labeling the branch points
  $a_s$ -- see section \ref{LabelBranchPoint}. For the case ${\cal B}$ it is important whether
  $[-1,1]$ intersects $[c_2,c_1]$ or not. So we consider two subcases:
\be
\begin{array}{cl}
Case~~{\cal B}1: & [-1,1]\cap[c_2,c_1]=\emptyset,\\
Case~~{\cal B}2: & [-1,1]\cap[c_2,c_1]\neq\emptyset.\\
\end{array}
\ee
The case ${\cal B}2$ in turn is subdivided into three subcases:
\be
\begin{array}{cl}
Case~~{\cal B}21: & [-1,1]\subset[c_2,c_1],\\
Case~~{\cal B}22: & [-1,1]\supset[c_2,c_1],\\
Case~~{\cal B}23: & all~ the~ rest.\\
\end{array}
\ee

\subsection{Pair of Pants associated to $R_3$}
\label{PantsOfR3}
  For the obvious reason, {\it a pair of pants} is the name for the
  Riemann sphere with three holes in it. Pair of pants may
  be conformally mapped to $\hat{\mbox{\sets C}}$ with  three
  nonintersecting {\it real} slots. This mapping is unique up to real
  linear-fractional mappings. The conformal
  class of pants with labeled boundary components depends on three
  real parameters varying in a cell.

{\bf D~e~f~i~n~i~t~i~o~n:}
To every PS-3 equation we associate the pair of pants:
\be
\label{R3Pants}
{\cal P}(R_3):=Cl
\left(
\hat{\mbox{\sets C}}
\setminus
\{
([-1,1]\triangle [a_1,a_2])\cup [a_3,a_4]
\}
\right)
\ee
  where $\triangle$ is the symmetric difference (union of two
  sets minus their intersection). {\it Closure} here and everywhere below
  is taken with respect to the intrinsic spherical metric
  when every slot acquires two banks. Boundary components of
  pants are colored in accordance with the palette:
\begin{center}
\begin{tabular}{ll}
$[-1,1]\setminus [a_1,a_2]$& -- red\\
$[a_1,a_2]\setminus [-1,1]$& -- blue\\
$[a_3,a_4]$& -- green
\end{tabular}
\end{center}

Thus obtained pair of pants will have boundary ovals of all three colors,
but in cases ${\cal B}21$ (green and two blue ovals) and ${\cal B}22$
(green and two red ovals). In case ${\cal A}$ the red, green and blue slots
always follow in the natural cyclic order of the extended real axis.

\subsection{Gauge Transformations}
Let us recall that parameter $R(x)$ is not uniquely determined by two
domains $\Omega_1$ and $\Omega_2$. Composition with linear-fractional transformations
preserving the segment $[-1,1]$ is admissible. The general appearance of such a mapping is

\begin{equation}
L^\pm_\alpha(t):=\pm\frac{t+\alpha}{\alpha t+1},
\quad \alpha\in (-1,1).
\label{LFI}
\end{equation}

\begin{lmm}
1. The gauge transformation $R\to L_\alpha^\pm\circ R$
does not change neither eigenvalues $\lambda$ nor the eigenfunctions $u(t)$
of any PS integral equation.

2. The gauge transformation $R\to R\circ L_\alpha^\pm$
does not change the eigenvalues $\lambda$ and slightly changes the eigenfunctions: $u(t)\to u(L_\alpha^\pm(t))$.
\end{lmm}
P~r~o~o~f. To simplify the notations we put $L(t):=L_\alpha^\pm(t)$.

1. The gauge transformation just adds a constant term to the right hand side of equation.
$$
\int_{-1}^1
\frac{u(t)dL(R(t))}{L(R(t))-L(R(x))}=
\int_{-1}^1
\frac{u(t)L'(R(t))dR(t)}{[L'(R(t))L'(R(x))]^{1/2}(R(t)-R(x))}=
$$
$$
\int_{-1}^1
\frac{u(t)dR(t)}{R(t)-R(x)}-
\int_{-1}^1
\frac{u(t)dR(t)}{R(t)-L^{-1}(\infty)}.
$$

2. We define the new variable $x_*:=L(x)$ and the new function $u_*(x_*):=u(x)$.
$$
\int_{-1}^1
\frac{u(t)dR(L(t))}{R(L(t))-R(L(x))}=
\pm\int_{-1}^1
\frac{u_*(t_*)dR(t_*)}{R(t_*)-R(x_*)},
$$
$$
\int_{-1}^1
\frac{u(t)dt}{t-x}=
\pm\int_{-1}^1
\frac{u_*(t_*)dL^{-1}(t_*)}{L^{-1}(t_*)-L^{-1}(x_*)}=
\pm\int_{-1}^1
\frac{u_*(t_*)dt_*}{t_*-x_*}
\mp\int_{-1}^1
\frac{u_*(t)dt}{t-L(\infty)}. ~~\bl
$$

We see that essentially the space of $PS-3$ equations has real
  dimension 3, the same as the moduli space of pants. It is easy to check the following:

\begin{itemize}
\item
Any gauge transformation of the parameter $R_3(x)$ does
not change the type (${\cal A}$, ${\cal B}1$, \dots) of integral equation.
\item
The transformation $R_3\to R_3\circ L_\alpha^\pm$ does not change the associated
pants and preserves the colors of the boundary ovals.
\item
Associated to functional parameter $L_\alpha^\pm\circ R_3$ are the pants
$L_\alpha^\pm{\cal P}(R_3)$. The colors of the boundary ovals are transferred
by  $L_\alpha^\pm$, but in one case. When the type of integral equation is $\cal A$,
the transformation $L_\alpha^-$ interchanges blue and green colours on the boundaries.
\end{itemize}

\subsection{Reconstruction of $R_3(x)$ from the Pants}
\label{R3fromPants}
The parameter $R_3(x)$ of integral equation may be reconstructed,
given the pants ${\cal P}(R_3)$ and the type ${\cal A}$, ${\cal B}1\dots$ of the equation.
One has to follow the routine described below.

\paragraph{Restore the Labeling of the Branch Points.} In case ${\cal B}2$
we temporarily paint the real segment separating two non-green slots
in blue. The (extended) blue segment is set to be $[a_1,a_2]$; the
green segment is $[a_3,a_4]$. The relabeling $a_1\leftrightarrow
a_2$ and $a_3\leftrightarrow a_4$ is eliminated by the natural
cyclic order of the points $a_1,a_2,a_3,a_4$ on the extended real
axis.

\begin{figure}
\centerline{\psfig{figure=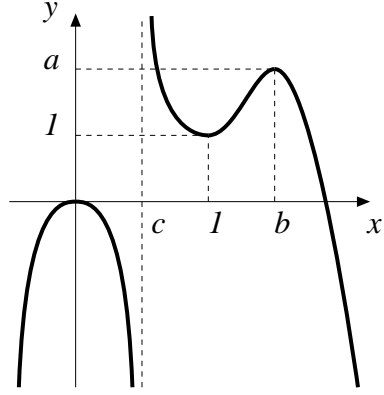}}
\caption[Two graphs]
{\small The graph of $\tilde {R}_3(x)$ when $c\in(\frac13,\frac12)$.}
\label{R3pm}
\end{figure}

 \paragraph{Normalized Covering.} Let $L_a$ be the unique linear-fractional map sending
  the points $a_1$, $a_2$, $a_3$, $a_4$ to respectively $0$,
  $1$, $a>1$, $\infty$. The conformal motion $L_b$ of the
  covering Riemann sphere sends the critical points $b_1$, $b_2$,
  $b_3$, $b_4$ of $R_3(x)$ (unknown at the moment) to respectively $0$, $1$,
  $b>1$, $\infty$. The function $L_a\circ R_3\circ
  L_b^{-1}$ with the normalized critical points and critical values
  takes a simple form:

  $$
  \widetilde{R_3}(x)=x^2L(x),
  $$
  with real linear fractional function $L(x)$ satisfying the
  restrictions:
$$
\begin{array}{ll}
L(1)=1,&L'(1)=-2,\\
L(b)=a/b^2,&L'(b)=-2a/b^3.
\end{array}
$$
We got four equations for three parameters of $L(x)$ and the unknown $b$.
The first two equations suggest the following expression for the
linear-fractional function:
$$
L(x)=1+2\frac{(c-1)(x-1)}{x-c}.
$$
The other two equations are solved parametrically in terms of $c$:
$$
b(c)=c\frac{3c-2}{2c-1};
\qquad
a(c)=c\frac{(3c-2)^3}{2c-1}.
$$
Given $a>1$, there are exactly two real solutions of the equation
$a(c)=a$. One of the solutions $c$ lies in the segment $(1/3,1/2)$, the other
lies in $(1,\infty)$. In the case $c\in(1/3,1/2)$ the segments $L_a(a_1,a_2)=(0,1)$
and $L_a(a_3,a_3)=(a,\infty)$ are filled with the points of the type (1:2),
which corresponds to our choice of labelling the branchpoints in section \ref{TopBrCov}.
The solution $c\in(1,\infty)$ is a fake as the same segments bear the points
of the type (3:0). Both functions $b(c)$ and $a(c)$ increase from $1$ to
$\infty$ when the argument $c$ runs from $1/3$ to $1/2$.

\paragraph{Reconstruction of the Mapping $L_b$.}
    In the case ${\cal A}$ the red, green and blue slots follow in the natural
    cyclic order. Hence, the segment $L_a[-1,1]$ is a subset of $[1,a]$. We
    choose the unique component of the pre-image $\widetilde{R_3}^{-1}$ of the
    segment $L_a[-1,1]$ belonging to the segment  $[1,b]$ -- see Fig.
    \ref{R3pm}. For the case ${\cal B}$ the segment $L_a[-1,1]$ is a subset of
    $(-\infty,a]$ and we choose  the pre-image of this segment which lie in
    $[b,\infty]$. The requirement: $L_b$ {\it maps $[-1,1]$ to the chosen
    segment} determines this map up to a pre-composition with the function
    (\ref{LFI}).\\[1mm]

We see that given the pants ${\cal P}(R_3)$, the functional parameter is
recovered up to a gauge transformation $R_3\to R_3\circ L_\alpha^\pm$. It is not
difficult to check, that the described above procedure applied to the pair of
pants $L_\beta^\pm{\cal P}(R_3)$ (in case $\cal A$ and the mapping $L_\beta^-$
reversing the orientation of real axis we additionally have to exchange the blue
and green colors of the slots) returns the covering map $L_\beta^\pm\circ
R_3\circ L_\alpha^\pm$. Roughly speaking, the classes of gauge transformation
of $R_3(x)$ correspond to the conformal classes of pants with suitably colored
boundary components and each conformal class of pants corresponds to a class of
certain functional parameter $R_3(x)$.

\section{Reduction to Projective Structures}
PS integral equations possess rich geometrical content
\cite{Bog1,Bog2} which we disclose in this section.
We describe a three-step reduction
  of the integral equation to a certain problem \cite{Bog1} for  projective
  structures on a riemann surface which has essentially combinatorial solution.

\subsection{Step 1: Functional Equation}

  Let us expand the kernel of the second integral in (\ref{PSE})
  into a sum of elementary fractions:

\begin{equation}
\label{kernel}
\frac{ R_3'(t)}{R_3(t)-R_3(x)}=
\frac{d}{dt}\log(R_3(t)-R_3(x))=
\sum_{k=1}^3\frac{1}{t-z_k(x)}-\frac{Q'}{Q}(t),
\end{equation}
where $Q(t)$ is the denominator in noncancellable representation of $R(t)$
as the ratio of two polynomials; $z_1(x)=x,~z_2(x),~z_3(x)$ --
are all solutions (including multiple and infinite) of the equation
$R_3(z)=R_3(x)$. This expansion suggests to rewrite the original equation (\ref{PSE})
as certain relationship for the Cauchy-type integral

\begin{equation}
\Phi(x):=\int\limits_I\frac{u(t)}{t-x}dt+const^*,
\quad x\in\hat{\mbox{\sets C}}\setminus[-1,1].
\label{Cauchy}
\end{equation}
Known $\Phi(x)$, the solution $u(t)$ may be recovered by the Sokhotskii-Plemelj
formula:
\begin{equation}
\label{SP}
u(t)=(2\pi i)^{-1}
\left[
\Phi(t+i0)-\Phi(t-i0)
\right],
\quad t\in I.
\end{equation}
The constant term $const^*$ in (\ref{Cauchy}), which we assume to be
\begin{equation}
const^*:=\frac{1}{\lambda - 3}
\left[
\int\limits_I\frac{u(t)
Q'(t)}{Q(t)}dt - const
\right]
\label{Const}
\end{equation}
is introduced to cancel the constant terms arising after
substitution of expression (\ref{Cauchy}) to the equation (\ref{PSE}).
In this way the following result was proven \cite{Bog2}:

\begin{lmm}
For $\lambda\neq 1,3$ the transformations (\ref{Cauchy}) and
(\ref{SP}) bring about a one-to-one correspondence between the
H\"older eigenfunctions $u(t)$ of PS-3 integral equation and the
holomorphic on the Riemann sphere outside the slot $[-1,1]$
nontrivial solutions $\Phi(x)$ of the functional equation
\begin{equation}
\label{FE}
\Phi(x+i0)+\Phi(x-i0)=\delta
\biggl(
\Phi (z_2(x))+\Phi(z_3(x))
\biggl)~,
\quad x\in I,
\end{equation}

\begin{equation}
\label{delta}
\delta =2/(\lambda -1),
\end{equation}
with H\"older boundary values $\Phi(x\pm i0)$ at the banks of the slot $[-1,1]$.
\end{lmm}

\subsection{Step 2: Riemann Monodromy Problem}
\label{RHP}
  In this section we reduce our functional (and
  therefore integral) equation to the Riemann monodromy
  problem in the following form. {\it Find a holomorphic vector
  $W(y)=(W_1,W_2,W_3)^t$ on the slit sphere
  ${\cal P}(R_3)\setminus[-1,1]$
  whose boundary values on the opposite sides of
  every slot are related by the constant matrix specified for each slot.}

\subsubsection{Monodromy Generators}
To formulate the Riemann monodromy problem we introduce $3\times 3$ permutation matrices
\be
\label{permutations}
{\bf D_1}:=
\begin{array}{||lll||}
1&0&0\\0&0&1\\0&1&0\\
\end{array}~;
\qquad
{\bf D_2}:=
\begin{array}{||lll||}
0&0&1\\0&1&0\\1&0&0\\
\end{array}~;
\qquad
{\bf D_3}:=
\begin{array}{||lll||}
0&1&0\\1&0&0\\0&0&1\\
\end{array}
\qquad
\ee
and a matrix depending on the spectral parameter $\lambda$:
\begin{equation}
\label{D}
{\bf D}:=
\begin{array}{||ccc||}
-1       & \delta & \delta \\
 0       & 1      & 0  \\
 0       & 0      & 1  \\
\end{array}
~,
\qquad 
\delta =2/(\lambda -1).
\end{equation}

\begin{lmm}
  Matrices ${\bf D_1}$, ${\bf D_2}$, ${\bf D_3}$,
  ${\bf D}$, ${\bf D_1D}=$ ${\bf DD_1}$ have order
  two as $GL_3$ group elements.
  \label{DsProperty}
\end{lmm}

\subsubsection{Separating Branches of $R_3^{-1}$}

  Let domain $\cal O$ be the compliment to the segments
  $[a_1,a_2]$ and $[a_3,a_4]$ on the Riemann sphere.
  The pre-image $R_3^{-1}{\cal O}$
  consists of three components ${\cal O}_j$, $j=1,2,3$, mapped one-one to
  $\cal O$ -- see Fig. \ref{CoverTopology}.
  Two of the components are (topological) discs with
  a slot and the third is an annulus. The enumeration of domains
  ${\cal O}_j$ is determined by the following rule:
  the segment $[-1,1]$ lies in the closure of ${\cal O}_1$,
the segment $[c_3,c_4]$ lies on the border of ${\cal O}_2$.

\subsubsection{}
  Let $u(x)$ be the solution of integral equation (\ref{PSE}) in the case ${\cal A}$.
  We consider the vector
\be
\label{Wofy}
W(y)=(\Phi(x_1),\Phi(x_2),\Phi(x_3))^t,
\qquad y\in{\cal O}\setminus[-1,1],
\ee
  where $\Phi(x)$ is from (\ref{Cauchy}) and $x_s$ is the unique
  solution of the equation $R_3(x_s)=y$ in ${\cal O}_s$. Vector
  $W(y)$ will be holomorphic and bounded in the domain ${\cal
  O}\setminus[-1,1]$  as all three points $x_s$, $s=1,2,3$,
  remain in the holomorphy domain of the function $\Phi(x)$. We claim
  that

\begin{center}
  \begin{tabular}{cc}
  $W(y+i0)={\bf D}W(y-i0)$,& when $y\in [-1,1]$,\\
  $W(y+i0)={\bf D_3}W(y-i0)$,& when $y\in [a_1,a_2],$\\
  $W(y+i0)={\bf D_2}W(y-i0)$,& when $y\in [a_3,a_4].$
  \end{tabular}
\end{center}

  Indeed, let $y^+:=y+i0$ and $y^-:=y-i0$ be two points on the opposite
  banks of $[a_1,a_2]$. Their inverse images $x_3^+=x_3^-$, $x_1^\pm=x_2^\mp$
  lie outside the cut $[-1,1]$. Hence $W(y^+)=$
  ${\bf D_3}W(y^-)$. For two points $y^\pm$ lying on the opposite
  banks of the slot $[a_3,a_4]$, their inverse images satisfy
  the relations $x_2^+=x_2^-$, $x_1^\pm=x_3^\mp$, which means
  $W(y^+)=$ ${\bf D_2}W(y^-)$. Finally, let $y^\pm$ lie on the  banks
  of $[-1,1]$. Now two points $x_2^+=x_2^-$ and $x_3^+=x_3^-$ lie
  in the holomorphy domain of $\Phi(x)$ while $x_1^+$ and $x_1^-$
  appear on the opposite sides of the cut $[-1,1]$. According to
  the functional equation (\ref{FE}),
$$
\Phi(x_1^+)=
-\Phi(x_1^-)+\delta(\Phi (x_2^-)+\Phi (x_3^-)),
$$
therefore $W(y^+)=$ ${\bf D}W(y^-)$ holds on the slot $[-1,1]$.

\subsubsection{}
  Conversely, let $W(y)=(W_1,W_2,W_3)^t$ be the bounded solution of the Riemann monodromy
  problem stated above. We define a piecewise holomorphic
  function  on the Riemann sphere:
  \be
  \Phi(x):=W_s(R_3(x)),
  \qquad {\rm when}~~x\in {\cal O}_s\setminus R_3^{-1}[-1,1],
  \qquad s=1,2,3.
  \label{PhifromW}
  \ee
  From the boundary relations for the vector $W(y)$ it
  immediately follows that the function $\Phi(x)$ has no jumps on
  the lifted cuts $[a_1,a_2]$, $[a_3,a_4]$, $[-1,1]$
  apart from the cut $[-1,1]$ from the upper sphere. Say, if the two points
  $y^\pm$ lie on the opposite sides of the cut $[a_1,a_2]$, then
  $W_3(y^+)=W_3(y^-)$ and  $W_1(y^\pm)=W_2(y^\mp)$ which means
  that the function $\Phi(x)$ has no jump on the components of
  $R_3^{-1}[a_1,a_2]$. From the boundary relation on the cut $[-
  1,1]$ it follows that $\Phi(x)$ is the solution for the functional
  equation (\ref{FE}). Therefore it gives a solution of Poincare-
  Steklov integral equation with parameter $R_3(x)$. Combining
  formulae (\ref{SP}) with (\ref{PhifromW}) we get the
  reconstruction rule
\be
\label{uofx}
u(x)=(2\pi i)^{-1}\biggl(W_1(R_3(x)+i0)-W_1(R_3(x)-i0)\biggr),
\qquad x\in [-1,1].
\ee

  \subsubsection{}
We have just proved for the case ${\cal A}$ the following
\begin{thrm} \cite{Bog1}\label{PSeqRMP}
~~If $\lambda\neq1,3$ then two formulas (\ref{Wofy}) and (\ref{uofx})
implement the one-to-one  correspondence between the solutions $u(x)$ of
the integral equation (\ref{PSE}) and the bounded solutions $W(y)$ of the
Riemann monodromy problem in the slit sphere $\hat{\mbox{\sets C}}\setminus\{
[a_1,a_2]\cup[a_3,a_4]\cup[-1,1]\}$ with the following matrices assigned to the slots:\\
\be
\label{Rbvp}
\begin{array}{c|c|c|c|c|}
&[-1,1]\setminus[a_1,a_2]&[a_1,a_2]\setminus[-1,1]&[a_3,a_4]&[-1,1]\cap[a_1,a_2]\\
\hline
{\rm Case~{\cal A}:}&{\bf D}& {\bf D_3}& {\bf D_2}& \\
{\rm Case~{\cal B}1:}&{\bf D}& {\bf D_1}& {\bf D_2}& \\
{\rm Case~{\cal B}2:}&{\bf D}& {\bf D_1}& {\bf D_2}& {\bf D_1D}={\bf DD_1}\\
\end{array}
\ee
\end{thrm}

\subsubsection{Monodromy Invariant}
It may be checked that the matrices $\bf D$, ${\bf D_1}$, ${\bf D_2}$,
${\bf D_3}$ generating the monodromy group for the solution $W(y)$
 are pseudo-orthogonal, that is preserve the same quadratic form
 \be
\label{J}
J(W):=\sum\limits_{k=1}^3 W_k^2-
\delta\sum\limits_{j<s}^3 W_jW_s.
\ee

  This form is not degenerate unless $-2\ne\delta\ne 1$, or
  equivalently $0\ne\lambda\ne 3$. Since the solution $W(y)$ of
  our monodromy problem is bounded near the cuts, the
  value of the form $J(W)$ is independent of the variable $y$.
  Therefore the solution ranges either in the smooth
  quadric $\{W\in\mbox{\sets C}^3:~~J(W)=J_0\neq0\}$, or the cone
  $\{W\in\mbox{\sets C}^3:~~J(W)=0\}$.

\subsubsection{Geometry of Quadric Surface}
\begin{figure}
\centerline{\psfig{figure=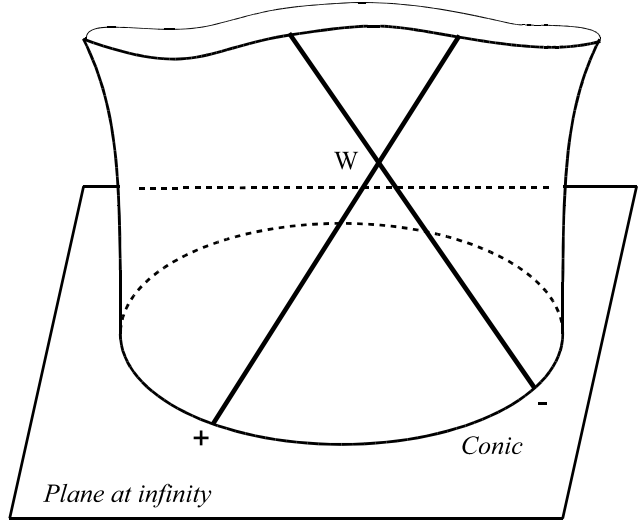}}
\caption[Quadric]
{\small Global coordinates $p^+$ and $p^-$ on quadric }
\label{Quadric}
\end{figure}

The nondegenerate projective quadric $\{J(W)=J_0\}$ contains two families of
line elements which for convenience we denote by the signs $'+'$ and $'-'$.
Two different lines from the same family are disjoint while two lines from
different families  intersect. The intersection of those lines with
the 'infinitely distant' secant plane gives points on the conic
\be
\label{conica}
\{(W_1:W_2:W_3)^t\in\mbox{\sets CP}^2:\quad J(W)=0\}
\ee
which by means of stereographic projection $p$ may be identified with the Riemann sphere.
Therefore we have introduced two global coordinates $p^\pm(W)$ on the quadric,
'infinite part' of which (i.e. conic (\ref{conica})) corresponds to coinciding coordinates:
$p^+=p^-$ (see Fig. \ref{Quadric}).

To obtain explicit expressions for the coordinate change $W\leftrightarrow p^\pm$ on quadric
we bring the quadratic form $J(W)$ to the simpler form
$J_\bullet(V):=$ $V_1V_3-V_2^2$ by means of the linear coordinate change
\be
W={\bf K}V
\label{V}
\ee
where
\be
\label{K}
{\bf K}:=
(3\delta+6)^{-1/2}
\begin{array}{||lll||}
1&1&1\\
1&\varepsilon^2&\varepsilon\\
1&\varepsilon&\varepsilon^2
\end{array}
\cdot
\begin{array}{||lll||}
0&\mu^{-1}&0\\
0&0&1\\
1&0&0
\end{array}
~,
\ee
$$
\varepsilon:=\exp(2\pi i/3),
\qquad
\mu:=\sqrt{\frac{\delta-1}{\delta+2}}=
\sqrt{\frac{3-\lambda}{2\lambda}}.
$$
Translating the first paragraph of the current section
into the language of formulae we get
\be
\label{V2Ppm}
p^\pm(W):=\frac{V_2\pm i\sqrt{J_0}}{V_1}
=\frac{V_3}{V_2\mp i\sqrt{J_0}};
\ee
and inverting this dependence,
\be
\label{Ppm2W}
W(p^+,p^-)=\frac{2i\sqrt{J_0}}{p^+-p^-}{\bf K}
\left(
\begin{array}{c}
1\\
(p^++p^-)/2\\
p^+p^-
\end{array}
\right)~.
\ee
The point $W(p^+,p^-)$ with coordinate $p^+$ (resp. $p^-$) being fixed moves along the straight
line with the directing vector ${\bf K}(1:p^+:(p^+)^2)$ (resp. ${\bf K}(1:p^-:(p^-)^2)$)
belonging to the conic (\ref{conica}).

\begin{lmm}
There exists a (spinor) representation $\chi:\quad O_3(J)\to PSL_2(\mbox{\sets C})$ such that:\\
1) The restriction of $\chi(\cdot)$ to $SO_3(J)$ is an isomorphism to $PSL_2(\mbox{\sets C})$.\\
2) For coordinates $p^\pm$ on the quadric the following transformation rule holds:\\
\be
\label{Ptrans}
\begin{array}{ll}
p^\pm({\bf T}W)=\chi({\bf T})p^\pm(W),&
\qquad{\bf T}\in SO_3(J),\\
p^\pm({\bf T}W)=\chi({\bf T})p^\mp(W),&
\qquad{\bf T}\in O_3(J)\setminus SO_3(J).\\
\end{array}
\ee
3) The linear-fractional mapping $\chi p:= (ap+b)/(cp+d)$ is the image of the matrix:
\be
\label{chiT}
{\bf T}:=\frac1{ad-bc}~{\bf K}~
\begin{array}{||ccc||}
d^2&2cd&c^2\\
bd&ad+bc &ac\\
b^2&2ab&a^2
\end{array}
~{\bf K}^{-1}
\in SO_3(J).
\ee
4) The generators of the monodromy group are mapped to the following elements:
\be
\label{chiD}
\begin{array}{ll}
\chi({\bf D}_s)p=\varepsilon^{1-s}/p,&\qquad s=1,2,3;\\
\chi({\bf D})p=
\displaystyle{\frac{\mu p-1}{p-\mu}}.&
\end{array}
\ee
\end{lmm}
P~r~o~o~f: We define the action of matrix ${\bf A}\in
SL_2(\mbox{\sets C})$ on the vector $V\in\mbox{\sets C}^3$ by the
formula:

\be
\label{SL2action}
{\bf A}:=\begin{array}{||cc||}a&b\\c&d\end{array}:
\qquad
\begin{array}{||cc||}
V_3&V_2\\
V_2&V_1
\end{array}
\quad
\longrightarrow
\quad
{\bf A}~
\begin{array}{||cc||}
V_3&V_2\\
V_2&V_1
\end{array}
~{\bf A}^t.
\ee
  It is easy to check that (\ref{SL2action}) gives the faithful
  representation of connected 3-dimensional group
  $PSL_2(\mbox{\sets C}):=SL_2(\mbox{\sets C})/\{\pm {\bf 1}\}$ into
  $SO_3(J_\bullet)$ and therefore, an isomorphism. Let us  denote
  $\chi_\bullet$ the inverse isomorphism  $SO_3(J_\bullet)\to$
  $PSL_2(\mbox{\sets C})$ and let $\chi(\pm {\bf
  T}):=\chi_\bullet({\bf K}^{-1}{\bf TK})$ for ${\bf T}\in
  SO_3(J)$. The obtained homomorphism $\chi:~~$ $O_3(J)\to$
  $PSL_2(\mbox{\sets C})$ will satisfy statement 1) of the lemma.

  To prove 2) we replace vector $V$ components in the right-hand
  side of (\ref{SL2action}) with their representation in terms of
  the stereographic coordinates $p^\pm$:

$$
\frac{i\sqrt{J_0}}{p^+-p^-}
{\bf A}~
\left[
(p^+,1)^t\cdot(p^-,1)+
(p^-,1)^t\cdot(p^+,1)
\right]
~{\bf A}^t=
$$
$$
i\sqrt{J_0}
\frac{(cp^++d)(cp^-+d)}{p^+-p^-}
\left[
(\chi p^+,1)^t\cdot(\chi p^-,1)+
(\chi p^-,1)^t\cdot(\chi p^+,1)
\right]
=
$$
$$
\frac{i\sqrt{J_0}}{\chi p^+-\chi p^-}
\left[
(\chi p^+,1)^t\cdot(\chi p^-,1)+
(\chi p^-,1)^t\cdot(\chi p^+,1)
\right]
=
$$
$$
\begin{array}{||cc||}
V_3(\chi p^+,\chi p^-)&V_2(\chi p^+,\chi p^-)\\
V_2(\chi p^+,\chi p^-)&V_1(\chi p^+,\chi p^-)
\end{array}
~~,
$$
  where we set $\chi p:=(ap+b)/(cp+d)$. Now (\ref{Ptrans})
  follows immediately for ${\bf T}\in SO_3(J)$. It remains to
  check the transformation rule for any matrix ${\bf T}$ from
  the other component of the group $O_3(J)$, say ${\bf T}=-{\bf 1}$.

  Writing the action (\ref{SL2action}) component-wise we arrive
  at conclusion 3) of the lemma.

  An finally, expressions 4) for the generators of monodromy
  group may be obtained either from analyzing formula
  (\ref{chiT}) or finding the eigenvectors of the matrices ${\bf
  D}_s,{\bf D}$ which correspond to the fixed points of linear-
  fractional transformations. ~~~\bl

\subsection{Step 3: Projective Structures}
Speaking informally, complex projective structure
  \cite{G,Tu,GKM,Hej,Man} on the Riemann surface $M$ is a
  meromorphic function $p$ on the universal cover $\widetilde{M}$ of
  the surface that transforms linear fractionally under the cover
  transformations. The appropriate representation $\pi_1(M)\to $
  $PSL_2(\mbox{\sets C})$ is called the {\it monodromy} of the
  structure $p$. The projective structure is called branched when
  $p$ has critical points. The set of
  all critical points of $p(t)$ with their multiplicities survives
  under the cover transformations of $\widetilde{\cal M}$.
  The projection of this set to the Riemann surface ${\cal M}$ is known as
  the  {\it branching divisor} ${\sf D}(p)$ of projective structure
  and the branching number of the structure $p(t)$ is $\deg{\sf D}(p)$.
  The classical examples of unbranched projective structures
  arise in Fuchsian or Schottky uniformization of Riemann
  surfaces.

\subsubsection{Projective Structures Generated by Eigenfunction}
Stereographic coordinates $p^\pm(y):=p^\pm(W(y))$ for the solution
of the Riemann monodromy problem (\ref{Rbvp}) will give two nowhere
coinciding  meromorphic functions in the sphere with three possibly
overlapping slots. As it follows from the transformation rules
(\ref{Ptrans}), the boundary values of two functions $p^\pm(y)$ on
the slot $D_*$, one of $[a_1,a_2]\setminus [-1,~1]$,
$[-1,~1]\setminus [a_1,a_2]$, $[a_1,a_2]\cap[-1,1]$ or $[a_3,a_4]$,
are related by the formulas \be \label{BVp}
\begin{array}{ll}
p^\pm(y+i0)=\chi({\bf D}_*)p^\mp(y-i0),&
\qquad y\in D_*\neq[a_1,a_2]\cap[-1,~1],\\
p^\pm(y+i0)=\chi({\bf DD_1})p^\pm(y-i0),&
\qquad y\in [a_1,a_2]\cap[-1,~1],
\end{array}
\ee
where ${\bf D}_*$ is the matrix assigned to the slot $D_*$ in (\ref{Rbvp}).

Relations (\ref{BVp}) allow us to analytically continue both
functions $p^+(y)$ and $p^-(y)$ through any slot to locally single
valued functions on the genus 2 Riemann surface 
\be
M:=\{w^2=(y^2-1)\prod\limits_{s=1}^4(y-a_s)\}, 
\label{M} 
\ee \
since all matrices ${\bf D}_*$ are involutive -- see Lemma \ref{DsProperty}. 
Further continuation gives single valued functions $p^\pm(\cdot)$ on the universal 
covering $\tilde{M}$. Traveling of the argument $y$ along the handle of the 
surface $M$ may result in the linear- fractional transformation of the value 
$p^\pm(y)$. Say, the continuations of $p^+(y)$ from the pants through the red 
and green slots will give two different functions on the second sheet related by 
the linear-fractional mapping $\chi({\bf DD_2})$.

\subsubsection{Branching of Structures $p^\pm$}
The way we have carried out the continuation of functions $p^\pm(y)$ suggests
that the branching divisors of the arising projective structures are related via
formula:
\be
\label{BranchSymm}
{\sf D}(p^+)=H{\sf D}(p^-)
\ee
where  $H(y,w):=(y,-w)$ is the hyperelliptic involution of the surface $M$.
We determine the branching numbers of the structures in the proof of

\begin{thrm} \cite{Bog1} \label{J0neq0}
  When $\lambda\not\in\{0,1,3\}$ the solutions $u(x)$ of the PS-3 integral equation
  are in one-one correspondence with the
  couples of meromorphic in the  slit sphere
  $\hat{\mbox{\sets C}}\setminus\{[a_1,a_2]\cup[a_3,a_4]\cup[-1,1]\}$ functions
  $p^\pm(y)$ with boundary values satisfying (\ref{BVp}) and either non or two
  critical points in common. The correspondence $u(x)\to$ $p^\pm(y)$ is
  implemented by the sequence of formulae (\ref{Cauchy}), (\ref{Wofy}) and
  (\ref{V2Ppm}). The inverse dependence is given up to proportionality 
  by the formula
\be
\label{UniU}
u(x)=\sqrt{\frac{\Omega(y)}{dp^+(y)dp^-(y)}}
(p^+(y)p^-(y) -\mu(p^+(y)+p^-(y))+1),
\ee
  where $x\in[-1,1]$ and $y:=R_3(x)+i0$,
  $\displaystyle{\Omega(y)=(y-y_1)(y-y_2)\frac{(dy)^2}{w^2(y)}}$
  is the holomorphic quadratic differential on
  the Riemann surface $M$ with zeroes at the critical points of
  the (possibly coinciding) functions $p^+$ and $p^-$, or with
  double zeroes $y_1=y_2$ (otherwise arbitrary) when $p^+=p^-$ is unbranched.
\end{thrm}

{\bf Remark}:  The number of the critical points of the structures
in the slit sphere is counted with their {\it weight and
multiplicity}: ~~1) the branching number of $p^\pm(y)$ at the branch
point  $a\in \{\pm1, a_1,\dots,a_4\}$ of $M$ is computed with
respect to the local parameter $z=\sqrt{y-a}$, ~~2) every branch
point of projective structure on the boundary of the pants should be
considered as a half-point.

{P~r~o~o~f}: {\bf 1.~~} Let $u(x)$ be an eigenfunction of
integral equation PS-3. The stereographic coordinates $p^\pm(y)$ of
the solution $W(y)$ of the associated Riemann monodromy problem are nowhere equal 
meromorphic functions when the invariant $J_0\neq0$, or two identically equal
functions when $J_0=0$. In any case they inherit the
boundary relationship (\ref{BVp}). 

What remains is to find the
branching numbers of the entangled structures $p^\pm(y)$. To this
end we consider the $O_3(J)$-invariant quadratic differential form 
$J(dW)=J_\bullet(dV)$ transferred to the slit sphere.

In the general case $J_0\neq0$ we get (up to proportionality) the 
Kleinian quadratic differential:
\be 
\label{Klein} 
\Omega(y)=\frac{dp^+(y)dp^-(y)}
{(p^+(y)-p^-(y))^2}, \qquad y\in\hat{\mbox{\sets C}}. 
\ee 
This
expression is the infinitesimal form of the cross ratio, hence it
remains unchanged after the same linear-fractional transformations
of the functions $p^+$ and $p^-$. Therefore, (\ref{Klein}) is a well
defined quadratic differential on the entire sphere. Lifting
$\Omega(y)$ to the surface $M$ we get a holomorphic differential.
Indeed, $p^+\neq p^-$ everywhere and applying suitable
linear-fractional transformation we assume that
$p^+=1+z^{m_+}+\{terms~of~higher~order\}$ and $p^-=cz^{m_-}+...$ in
terms of local parameter $z$ of the surface, $m_\pm\ge1$, $c\neq0$.
Then $\Omega=cm_+m_-z^{m_++m_--2}(dz)^2+\{terms~of~higher~order\}$.
Therefore
$$
{\sf D}(p^+)+{\sf D}(p^-)=(\Omega).
$$
Any holomorphic quadratic differential on genus 2 surface has 4 zeroes. 
The curve $M$ consists of two copies of the slit sphere interchanged by the 
hyperelliptic involution $H$. Taking
into account the symmetry (\ref{BranchSymm}) of branching divisors,
we see that the structures $p^\pm$ together have two critical points in the 
slit sphere.

In the special case $J_0=0$ two structures merge: $p^\pm(y)=:p(y)$ and the 
same quadratic differential $J(dW)=J_\bullet(dV)$ on the curve $M$ has the appearance:
\be
\label{OmegaCone}
\Omega(y)= [V_1(y)dp(y)]^2 ,
\ee
here $V_1(y)$ is the first component in the vector $V(y)$ defined by formula (\ref{V}).
The analysis of this representation in local coordinates suggests that
\be
2{\sf D}(p)+2(W)=(\Omega),
\ee
where $(W)$ is the divisor of zeroes of the locally holomorphic (but globally multivalued)
on $M$ vector $W(y)$. To characterize $(W)$ we need the following lemma, which we prove 
at the end of the current section. 

\begin{lmm}
\label{ZeroDivW}
The vector $W(y)$ cannot have simple zeroes at the fixed points of the hyperelliptic involution
of $M$ when $J_0=0$ and $\lambda\neq0,3$.
\end{lmm}

The divisor $(W)$ is obviously invariant under the hyperelliptic involution $H$. 
From this Lemma it follows that either $(W)=0$ (therefore $\deg {\sf D}(p)=2$) 
or $(W)$ consists of two points interchanged by $H$ (therefore the structure $p$ 
is unbranched). In other words, $p(y)$ has the branching number 0 or 2 on 
the slit sphere and the quadratic differential $\Omega$ is a square of a holomorphic 
linear differential.

{\bf 2.~~}
Conversely, let $p^+(y)$ and $p^-(y)$ be two not identically equal meromorphic 
functions on the slit sphere, with boundary conditions (\ref{BVp}) and the total 
branching number either zero or two (see remark above). For the meromorphic 
quadratic differential (\ref{Klein}) on the Riemann surface $M$ we establish 
(using local coordinate on the surface) the inequality:
\be
\label{DivIneq}
{\sf D}(p^+)+{\sf D}(p^-)\ge(\Omega)
\ee
where the deviation from equality means that there is a point where $p^+=p^-$.
But the degree of the divisor on the left of (\ref{DivIneq}) is zero or four and 
 $\deg(\Omega)=4$. Therefore, $p^+\neq p^-$ everywhere (and the total branching 
 of this pair of functions in the slit sphere is two). 
 
 The holomorphic vector $W(p^+(y),p^-(y))$ in the slit sphere solves the 
 Riemann monodromy problem specified in theorem \ref{PSeqRMP}. We already know 
 how to convert the latter vector to the eigenfunction of integral equation PS-3. 
Careful computation gives the restoration formula 
\be
\label{UfromP}
2\pi u(x)=\sqrt{\frac{(\delta+2)J_0}3}
\frac{p^+(y)p^-(y) -\mu(p^+(y)+p^-(y))+1}
{p^+(y)-p^-(y)},
\ee
where $x\in[-1,1]$ and $y:=R_3(x)+i0$. Formula (\ref{UniU}) appears after 
the substitution of (\ref{Klein}) to the latter formula.

Finally, suppose that two functions $p^\pm(y)$ with necessary branching and 
boundary behaviour are identical. For the solution on the cone,
$V=V_1~(1,p,p^2)^t$ and the first component $V_1$ may be taken from  
(\ref{OmegaCone}). Therefore we consider the vector on the slit sphere:
\be
\label{UniW}
W(y):=\frac{(y-y_1)}{w(y)p'(y)}
{\bf K}(1, p(y), p^2(y))^t,
\ee
  where $y_1$ is the critical point of $p(y)$ or arbitrary point when $p(y)$ is unbranched.
  One immediately checks that it is holomorphic and solves the Riemann monodromy problem 
  specified in theorem \ref{PSeqRMP}. Now to find the corresponding eigenfunction is a routine task.
~~~\bl

P~r~o~o~f~~o~f~~L~e~m~m~a~ \ref{ZeroDivW}.
Let $z$ be local coordinate on $M$ in the vicinity of the fixed point $z=0$ 
of the hyperelliptic involution $z\to -z$. 
Boundary relationship of the vector $W$ on the slots implies the symmetry:
\be
\label{WSymm}
W(-z)={\bf D_*}W(z)
\ee
where ${\bf D_*}$ is one of the matrices ${\bf D_1}$, ${\bf D_2}$, ${\bf D_3}$ or 
${\bf D}$. The matrix ${\bf D_*}$ has eigenvalues $+1$, $+1$, $-1$ and the 
intersection of the cone $\{J(W)=0\}$ with the invariant plane corresponding to 
to the eigenvalue $+1$ contains eigenvector $W^+$ of the matrix. Suppose that 
$W(z)=W^-z+\{terms~of~higher~order\}$. 
Substituting the last ansatz to (\ref{WSymm}) we see that $W^-$ 
is the eigenvector of the matrix corresponding to eigenvalue $-1$. Let $J(\cdot,\cdot)$
be the bilinear form polar to quadratic form $J(\cdot)$, then 
$$
J(W^+,W^-)=
J({\bf D_*}W^+,{\bf D_*}W^-)=
-J(W^+,W^-)=0.
$$
Once vector $W(z)$ varies in the cone, $J(W^-)=0$. Now we see that the cone contains 
the entire plane generated by the vectors $W^+$ and $W^-$. Therefore the cone is degenerate 
which only happens when $\lambda=0,3$. ~~~\bl

\subsection{Mirror Symmetry of Solution}
\label{MirrorSym}
  In what follows we are looking for {\it real} solutions $u(x)$
  of the integral equation (\ref{PSE}). There is no loss of the
  generality. Indeed, the restrictions on the monodromy of projective
  structures \cite{Bog1} imply that the spectrum
  of any PS-3 integral equation belongs to the segment $[0,3]$.
  Now both  real and imaginary parts of any complex eigenfunction
  $u(x)$ are the solutions of the integral equation.

  Real solutions $u(x)$ of the integral equation correspond to
  the solutions of Riemann monodromy problem with mirror symmetry:
  $W(\bar{y})=\overline{W}(y)$. This symmetry for the vector $V(y):={\bf
  K}^{-1}W(y)$  takes the form
$
V(\bar{y})=
(\overline{V_3}(y)~{\rm sign}(\delta+2),
(\overline{V_2}(y)~{\rm sign}(\delta-1),
(\overline{V_1}(y)~{\rm sign}(\delta+2))
$.
  The values $\delta+2$ and $\delta-1$ have the same sign as it follows
  from the range of spectral parameter $\lambda\in[0,3]$. Therefore, real
  solutions are split into two classes depending on the sign
  of $(\delta+2)J_0$:
  $$
\begin{array}{ll}
Symmetric~~ ((\delta+ 2)J_0>0),& p^\pm(\bar{y})=1/\overline{p^\pm}(y)\\
Antisymmetric~~ ((\delta+ 2)J_0\le0),& p^\pm(\bar{y})=1/\overline{p^\mp}(y) \\
\end{array},
\qquad y\in {\cal P}(R_3)\setminus[-1,~1],
$$

  In the remaining part of the article we give explicit parametrization of all
  antisymmetric solutions for the integral PS-3 equations of
  the considered type -- when six points $\pm1, a_1,\dots,a_4$ are
  real and pairwise distinct.

Restricting ourselves to the search of antisymmetric solutions we
have to find only one function in the pants, say $p(y):=p^+(y)$ while the remaining function may be
recovered from the mirror antisymmetry:
\be
\label{Antisymm}
p^-(y)=1/\overline{p^+(\bar{y})}.
\ee
On the boundary components 
of the slit sphere this function obeys the rule:
$$
p^+(y\pm i0)=
\chi({\bf D}_*)p^-(y\mp i0)=
\chi({\bf D}_*{\bf D}_1)\overline{p^+(y\pm i0)},
\qquad y\in D_*\neq
[-1,1]\cap[a_1,a_2].
$$
Therefore

\centerline{
\begin{tabular}{cc}
$p\in\hat{\mbox{\sets R}}$,& when ${\bf D}_*={\bf D}_1$;\\
$p\in\varepsilon\hat{\mbox{\sets R}}$,& when ${\bf D}_*={\bf D}_2$;\\
$p\in\varepsilon^2\hat{\mbox{\sets R}}$,& when ${\bf D}_*={\bf D}_3$;\\
\end{tabular}}

\noindent
and finally when ${\bf D}_*={\bf D}$, the value of $p$ lies on the circle
\be
\label{circle}
C:=
\{~p\in\mbox{\sets C}:
\quad |p-\mu^{-1}|^2=\mu^{-2}-1~~   \},
\qquad \mu:=\sqrt{\frac{3-\lambda}{2\lambda}}
\ee

As an immediate consequence of this observation we give a universal restriction for the
spectrum of our integral equation
\begin{lmm}
 Antisymmetric eigenfunctions correspond to eigenvalues
$\lambda\in [1,3]$.
\label{SpLocus}
\end{lmm}
P~r~o~o~f: For the cases ${\cal A}$, ${\cal B}1$, ${\cal B}22$, ${\cal B}23$ the
slot $[-1,1]\setminus[a_1,a_2]$ is not empty and the boundary value of  $p(y)$
on this slot belongs to the circle $C$. This circle is an empty set for $\mu>1$, or equivalently
$\lambda\in(0,1)$. The proof for the remaining case ${\cal B}21$ requires special machinery and
 will be given in Sect. \ref{B2Proof}. ~~\bl\\[5mm]

The critical points of two functions $p^+(y)$ and $p^-(y)$ in the considered antisymmetric case
are complex conjugate as it follows from  (\ref{Antisymm}). Taking this fact into account we 
reformulate Theorem  \ref{J0neq0} for the antisymmetric solutions:

\begin{thrm}
When $\lambda\not\in\{0,1,3\}$, the antisymmetric solutions $u(x)$ of
integral equation {\rm PS-3} are in one-two correspondence with the meromorphic 
in the slit sphere $\hat{\mbox{\sets C}}\setminus\{[a_1,a_2]\cup[a_3,a_4]\cup[-1,1]\}$
functions $p(y)$ that have either none or one critical point
in the domain and the following values on the boundary components:\\[1mm]

\begin{tabular}{c||c|c|c}
$y\in$ & $[-1,1]\setminus[a_1,a_2]$ (red)& $[a_1,a_2]\setminus[-1,1]$ (blue)& $[a_3,a_4]$ (green)\\
\hline
$p(y\pm i0)\in$& $C$ &
\begin{tabular}{c}
$\varepsilon^2\hat{\mbox{\sets R}}$ (Case ${\cal A}$)\\
$\hat{\mbox{\sets R}}$  (Case ${\cal B}$)
\end{tabular}
& $\varepsilon\hat{\mbox{\sets R}}$\\
\end{tabular}

\vspace{1mm}
\noindent
In case ${\cal B}2$ the function $p(y)$ has the jump on the remaining part of the boundary:
\be
\label{Jump}
p(y+i0)=\chi({\bf DD}_1)p(y-i0),
\quad y\in [-1,1]\cap[a_1,a_2].
\ee
\label{SemiMainTh}
\end{thrm}

{\bf Remark~}
By 'one-two' correspondence we mean the following: 
given any function $p(y)$ satisfying the conditions of this theorem,
it's easy to check that its {\it antisymmetrization} $1/\overline{p(\bar{y})}$
also satisfies all the conditions. Therefore, we have a correspondence of an eigenfunction 
$u(x)$ to a couple:  function $p(y)$ and its antisymmetrization, only one of them being independent.

\section{Statement of the main result}
\label{Result}

  From the Sect. \ref{MirrorSym} it follows that every
  antisymmetric eigenfunction $u(x)$ of PS-3 integral equation
  induces a mapping $p(y)$ of the pants ${\cal P}(R_3)$ to a
  multivalent domain spread possibly with branching over the
  Riemann sphere. Such surface is known as Kleinian membrane 
  or \"Uberlagerungsfl\"ache and the complete list of
  them is given in this section.

\subsection{Tailoring the Pants}

We define pants ${\cal PQ}(\lambda, h_1, h_2|m_1, \dots)$ of several
fashions $\cal Q$ which parametrically depend on spectral parameter
$\lambda$, two other reals $h_1,~h_2$ and one or two integers
$m_1,\dots$. Each boundary oval of our pair of pants covers a circle
and acquires its color in the following way:

\begin{center}
\begin{tabular}{ll}
$C$&     -- red,\\
$\varepsilon\hat{\mbox{\sets R}}$ or
$\chi({\bf DD}_1)\varepsilon\hat{\mbox{\sets R}}$ & -- green,\\
$\hat{\mbox{\sets R}}$ or $\varepsilon^2\hat{\mbox{\sets R}}$& -- blue.\\
\end{tabular}
\end{center}

Any constructed pair of pants may be obtained from the "basic" pants ${\cal
PQ}(\lambda, h_1, h_2|\dots)$ with lowest possible integer parameters by a
surgery procedure known as "grafting" and introduced independently by B.Maskit,
D.Hejhal and D.Sullivan--W.Thurston in 1969-1983.

\subsubsection{Cases ${\cal A,~B}1$}

For real $\lambda\in(1,2)$ we consider (depending on $\lambda$) open annulus
$\alpha$ bounded by two circles: $C$ defined in (\ref{circle}) and
$\varepsilon\hat{\mbox{\sets R}}$. Another annulus bounded by $C$ and
$\varepsilon^2\hat{\mbox{\sets R}}$ we denote $\bar{\alpha}$. Note that for
the considered values of $\lambda$ the circle $C$ does not intersect the lines
$\varepsilon^{\pm1}\mbox{\sets R}$. The $m-$ sheeted unbranched covering of the
annuli, $m=1,2,\dots$, we denote as $m\cdot\alpha$ or  $m\cdot\bar{\alpha}$
correspondingly.

The pants of four fashions ${\cal PA}_1$, ${\cal PA}_2$, ${\cal PA}_3$, ${\cal PB}1$
are sewn together of the annuli we have introduced in the way specified in Tab. \ref{PantsAB1}.

\begin{table}[h!]
\centering
\begin{tabular}
{l|l|l}
Fashion of Pants  & Range of $h_1,~h_2$ and $m_1,~m_2$ & {\bf Definition}\\
\hline
${\cal PA}_1(\lambda,h_1,h_2|~m_1,m_2)$&
\begin{tabular}{l}
$h:=h_1+ih_2\in\alpha\cap\bar{\alpha}$, $|h|\ge1$;
\\
$m_1,m_2=1,2,\dots$
\end{tabular}&
\begin{tabular}{l}
$Cl\{(m_1\cdot\alpha)\setminus[\mu^{-1}, h]\}+$\\
$Cl\{(m_2\cdot\bar{\alpha})\setminus[\mu^{-1}, h]\}$
\end{tabular}\\
\hline
${\cal PA}_2(\lambda,h_1,h_2|~m_1,m_2)$&
\begin{tabular}{l}
$0<h_1<h_2$, $h_1h_2\ge1$;\\
$m_1=1,2,\dots$, $m_2=0,1,2,\dots$
\end{tabular}&
\begin{tabular}{l}
$Cl\{(m_1\cdot\alpha)\setminus-\varepsilon^2[h_1,h_2]\}+$\\
$Cl\{m_2\cdot\bar{\alpha}\}$
\end{tabular}\\
\hline
${\cal PA}_3(\lambda,h_1,h_2|~m_1,m_2)$&
\begin{tabular}{l}
$0<h_1<h_2$, $h_1h_2\ge1$; \\
$m_1=0,1,2,\dots$, $m_2=1,2,3,\dots$
\end{tabular}&
\begin{tabular}{l}
$Cl\{(m_2\cdot\bar{\alpha})\setminus-\varepsilon[h_1,h_2]\}+$\\
$Cl\{m_1\cdot{\alpha}\}$
\end{tabular}\\
\hline
${\cal PB}1(\lambda,h_1,h_2|~m)$&
\begin{tabular}{l}
$\mu^{-1}+\sqrt{\mu^{-2}-1} <h_1<h_2$;\\
$m=1,2,3,\dots$
\end{tabular}&
$Cl\{(m\cdot\alpha)\setminus[h_1,h_2]\}$\\
\hline
\end{tabular}
\caption{Three-parametric families of pairs of pants for the cases $\cal A$, ${\cal B}1$;
parameter $1<\lambda<2$}
\label{PantsAB1}
\end{table}

\vspace{2mm}
Sign '+' in the definitions of Tab. \ref{PantsAB1} means certain surgery explained below.

\begin{figure}
\centerline{\psfig{figure=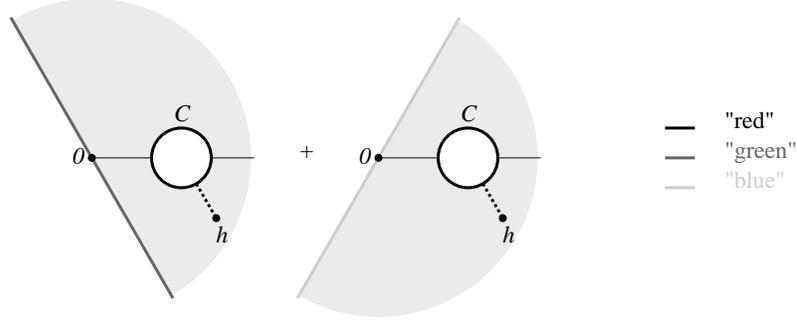}}
\caption[]
{\small The pair of pants  ${\cal PA}_1(\lambda,h_1,h_2|~m_1,m_2)$ is sewn down of two annuli}
\label{FA1}
\end{figure}

\begin{figure}
\centerline{\psfig{figure=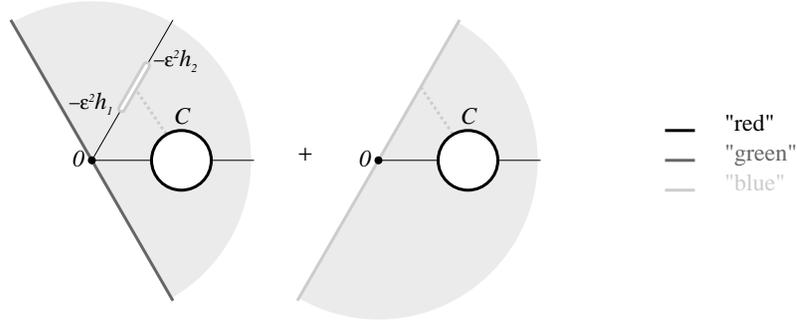}}
\caption[]
{\small The pair of pants ${\cal PA}_2(\lambda,h_1,h_2|~m_1,m_2)$ is sewn of
simpler pants and the annulus}
\label{FA2}
\end{figure}

\paragraph{Instructions on sewing annular patches together. }

{\bf 1.} ${\cal PA}_1(\lambda,h_1,h_2|~m_1,m_2)$.
  Take two annuli $m_1\cdot\alpha$ and  $m_2\cdot\overline{\alpha}$.  Cut the
  top sheet of each annulus along the same segment (dashed red line in the
  Fig. \ref{FA1}) starting at the point $h:=h_1+ih_2$ and ending at the circle
  $C$. Now sew the left bank of one cut on the right bank of the other. The
  emerging surface is the pair of pants.

{\bf 2.} ${\cal PA}_2(\lambda,h_1,h_2|~m_1,m_2)$.
    The annulus $m_1\cdot\alpha$ with the segment $-\varepsilon^2[h_1,h_2]$ removed
  from the top sheet is a pair of pants ${\cal PA}_2(\lambda,h_1,h_2|m_1,0)$. Cut
  the obtained pair of pants along the segment joining the circle $C$ to the
  slot (dashed blue line in the Fig. \ref{FA2}). Also cut top sheet of the annulus
  $m_2\cdot\overline{\alpha}$ along the same segment and sew the left bank of
  one cut on the right bank of the other. The arising surface is the
  pair of pants.

{\bf 3.} ${\cal PA}_3(\lambda,h_1,h_2|~m_1,m_2)$.
The annulus $m_2\cdot\bar{\alpha}$ with the segment $-\varepsilon[h_1,h_2]$ removed
from the top sheet, is a pair of pants ${\cal PA}_3(\lambda,h_1,h_2|0,m_2)$. As
in the previous passage, we may sew in the annulus $m_1\cdot{\alpha}$ to the
obtained pants to get the result.

{\bf 4.} ${\cal PB}1(\lambda,h_1,h_2|~m)$ is just the annulus $m\cdot\alpha$
with the segment $[h_1,h_2]$ removed  from its top sheet.

  The limit case of the pants ${\cal PA}_1$, when the branch point
  $h_1+ih_2$ tends to $\varepsilon^{\pm1}\mbox{\sets R}$, coincides with the
  limit cases of pants ${\cal PA}_2$ or ${\cal PA}_3$, when $h_1=h_2>0$.
The corresponding unstable two-parametric families of pants 
 ${\cal PA}_{12}$ and  ${\cal PA}_{13}$ are defined in Tab. \ref{InterMediate}.

\begin{table}[h!]
\centering
\begin{tabular}{l|l}
Fashion of Pants  & {\bf Definition}\\
\hline
${\cal PA}_{12}(\lambda,h|~m_1,m_2)$&
 \begin{tabular}{l}
  ${\cal PA}_1(\lambda, -Re(\varepsilon^2h),-Im(\varepsilon^2h)|m_1,m_2)=$\\
  ${\cal PA}_2(\lambda, h,h|m_1,m_2)$
 \end{tabular}\\
 \hline
${\cal PA}_{13}(\lambda,h|~m_1,m_2)$&
 \begin{tabular}{l}
  ${\cal PA}_1(\lambda, -Re(\varepsilon h),-Im(\varepsilon h)|m_1,m_2)=$\\
  ${\cal PA}_3(\lambda, h,h|m_1,m_2)$\\
  \end{tabular}
\end{tabular}
\caption{Unstable two-parametric families of pants.
The range of parameters: $1<\lambda<2$, $h>0$, $m_1$ and $m_2=1,2,3,\dots$.}
\label{InterMediate}
\end{table}

\subsubsection{Case ${\cal B}2$}
Two circles: $\varepsilon\hat{\mbox{\sets R}}$ and $\chi({\bf
DD}_1)\varepsilon\hat{\mbox{\sets R}}$ do not intersect when $\lambda\in(1,3)$.
They bound the open annulus $\beta$ depending on $\lambda$. The $m-$ sheeted
unbranched covering of the annulus  we denote as $m\cdot\beta$, $m=1,2,3,\dots$
The points of the annulus $m\cdot\beta$ may be described in the form
$$
p=\mu^{-1}+\rho\exp(i\phi),
$$
where $\rho>0$ and the argument $\phi\in$ {\sets R} $mod~2\pi m$.
The action of $\chi({\bf DD}_1)$ on the sphere (i.e. ~consecutive
reflections in circles $C$ and $\hat{\mbox{\sets R}}$) lifts to the
involution of the multi-sheeted annulus $m\cdot\beta$ in the
following way:
\be
\Xi: \mu^{-1}+\rho\exp(i\phi)\to\mu^{-1}+\frac{r^2}\rho\exp(-i\phi)
\label{mBetaInv} 
\ee
where $r:=\sqrt{\mu^{-2}-1}$ is the radius of the circle $C$.

{\bf Definition.} We introduce three pairs of pants ${\cal PB}21$, ${\cal PB}22$
and ${\cal PB}23$, each of them depend on three reals $\lambda$, $h_1$, $h_2$ and 
an integer $m$:
$$
{\cal PB}2s(\lambda,h_1,h_2|~m):=Cl\{
(m\cdot\beta)\setminus(E_s^1(h_1)\cup E_s^2(h_2))\}/\Xi,
\quad s=1,2,3,
$$
where the slots $E_s^1,E_s^2$ are defined in the Tab. \ref{SlotsDef}. The slots
are invariant with respect to the involution $\Xi$ and do not intersect each other as 
well as the boundary of the annulus $m\cdot\beta$.

\begin{table}[h!]
\centering
\begin{tabular}{l|l}
{\bf Definition} of slots & Range of $h_1,~h_2$\\
\hline
\begin{tabular}{l}
$E_1^1(h_1):=\mu^{-1}+r\exp[-h_1,h_1]$,\\
$E_1^2(h_2):=\mu^{-1}+r\exp[-h_2,h_2]\exp(i\pi m)$
\end{tabular}
&
\begin{tabular}{l}
$h_1\ge h_2>0$,\\ when $m$ is even;\\
$(\mu^{-1}+r\exp~h_1)\cdot$\\
$(\mu^{-1}-r\exp~h_2)\ge1$,\\ when $m$ is odd
\end{tabular}\\
\hline
\begin{tabular}{l}
$E_2^1(h_1):=\mu^{-1}+r\exp[-ih_1,ih_1]$,\\
$E_2^2(h_2):=\mu^{-1}+r\exp[-ih_2,ih_2]\exp(i\pi m)$
\end{tabular}
&
\begin{tabular}{l}
$h_1\ge h_2$, when $m$ is even;\\
${\rm Arg}(\exp(ih_1)+\mu r)\ge$\\
${\rm Arg}(\exp(ih_2)-\mu r)$\\
when $m$ is odd;\\
$h_1+h_2<m\pi$, $h_2>0$,\\
for any $m$
\end{tabular}\\
\hline
\begin{tabular}{l}
$E_3^1(h_1):=\mu^{-1}+r\exp[-h_1,h_1]$,\\
$E_3^2(h_2):=\mu^{-1}+r\exp[-ih_2,ih_2]\exp(i\pi m)$
\end{tabular}
& $h_1>0$, $m\pi>h_2>0$\\
\hline
\end{tabular}
\caption{Slots for the subcases of ${\cal B}2$, parameter $1<\lambda<3$.}
\label{SlotsDef}
\end{table}

To understand the interrelation of introduced constructions it is very useful 
to imagine how the pants ${\cal PB}1$
are transformed to the pants of fashion ${\cal Q=B}23$ and the latter --
to the pair of pants  ${\cal PB}21$ or ${\cal PB}22$.

\subsection{The main theorem}
\begin{thrm}

  When $\lambda\neq 1,3$ the antisymmetric eigenfunctions of {\rm PS-3} integral
  equation for the case ${\cal Q=A}$, ${\cal B}1$, ${\cal B}21$, ${\cal B}22$,
  ${\cal B}23$ are in one to one correspondence with the pairs of pants
  \footnote{For the case ${\cal A}$ there are three stable and two
  unstable pants fashions ${\cal PA}_*(\dots)$ } ${\cal
  PQ}(\lambda, h_1,h_2|m_1..)$ conformally equivalent to the  pants
  (\ref{R3Pants}) associated to the functional parameter of integral equation.

Let the function $p(y)$ conformally maps the pair of pants ${\cal P}(R_3)$  to the pants
${\cal PQ}(\lambda, h_1,h_2|m_1..)$, then up to proportionality

\begin{equation}
u(x)=
\left\{
\begin{array}{ll}
\displaystyle{
\sqrt{\frac{(y-y_1)(y-y_2)}{p'(y^+)p'(y^-)}}
\frac{p(y^+)-p(y^-)}{w(y)}},&\quad x\in[-1,1]\setminus[a_1,a_2],\\
\displaystyle{
\sqrt{(y-y_1)(y-y_2)}
\frac{{\rm Im}~p(y^+)}{w(y)|p'(y^+)|}},&\quad x\in[-1,1]\cap[a_1,a_2].
\end{array}
\right.
\label{UviaP}
\end{equation}

Here $y:=R_3(x)$, $y^\pm:=y\pm i0$. For the fashion ${\cal Q=A}_1$,
$y_1=\overline{y_2}$ is the inner critical point of the function $p(y)$;
for other fashions $\cal Q$ real $y_1$ and $y_2$ are boundary critical points
of the function $p(y)$.
\end{thrm}

The proof of the main theorem for the cases ${\cal A,B}1$ is given in Sect. \ref{AB1Proof}
and for the case ${\cal B}2$ -- in Sect. \ref{B2Proof}.

\subsection{Corollaries}
The representation (\ref{UviaP}) cannot be called explicit in the usual sense,
since it comprises a transcendent function $p(y)$. We show that nevertheless
the representation allows us to understand the following
properties of the solutions.

{\bf 1.} The "antisymmetric" part of the spectrum is always a subset of $[1,3]$;
for the equations of types ${\cal A,B}1$ this part of the spectrum always lies in
$[1,2]\cup\{3\}$.

{\bf 2.} Every $\lambda\in (1,3)$ is the eigenvalue for infinitely many equations PS-3.

P~r~o~o~f.~ Any of the constructed pants may be conformally mapped to the standard
form: the sphere with three real slots. Now we can apply the procedure of the Sect.
\ref{R3fromPants} and get the functional parameter $R_3(x)$ such that the
associated pair of pants is conformally equivalent to the pants we started from.

{\bf 3.} Eigenfunction $u(x)$ related to the pants ${\cal PQ}(\dots|m_1,m_2)$
has exactly $m_1+m_2+1$ zeroes on the segment $[-1,1]$ when ${\cal Q=A}$, ${\cal B}1$.

  P~r~o~o~f.~ According to the formula
  (\ref{UviaP}), the number of zeroes of eigenfunction  $u(x)$ is equal to
  the number of points $y\in[-1,1]$ where $p(y^+)=p(y^-)$. This
  number in turn is equal to the number of solutions of the
  inclusion

\begin{equation}
S(y):=
Arg [p(y^-)-\mu^{-1}] -
Arg [p(y^+)-\mu^{-1}]
\quad\in 2\pi\mbox{\sets Z},
\qquad y\in[-1,1].
\label{Inclus}
\end{equation}

  Let the point $p(y)$ goes $m$ times around the circle $C$  when
  its argument $y$ travels along the banks of $[-1,1]$. Integer $m$
  is naturally related to the integer parameters of pants ${\cal
  PQ}(\dots)$. The function $S(y)$ strictly increases from $0$
  to $2\pi m$ on the segment $[-1,1]$, therefore the inclusion
  (\ref{Inclus}) has exactly $m+1$ solutions on the mentioned
  segment. ~~~\bl

  {\bf 4.} The mechanism for arising the discrete spectrum of the integral
  equation is explained. Sewing annuli down to the pants ${\cal
  PQ}(\lambda,h_1,h_2|\dots)$ one changes the conformal structure of
  the latter. To return to the conformal structure specified by ${\cal P}(R_3)$ we have to
  change the real parameters of the pants, one of them is the
  spectral parameter $\lambda$.

  In a sense, the eigenvalue problem (\ref{PSE}) is reduced to
  the solution of three equations for three unknown numbers.
  These equations relate moduli of given pants ${\cal P}(R_3)$ to
  the moduli of membrane with real parameters $\lambda$, $h_1$,
  $h_2$ and extra discreet parameters.

\section{Auxiliary constructions}
\label{Construct}
The combinatorial analysis of the arising projective structure $p(y)$ is based on
two constructions we describe below.

  Let $p(y)$ be a holomorphic map from a Riemann surface $\cal M$
  with a boundary  to the sphere and the selected boundary
  component $(\partial {\cal M})_*$ is mapped to a circle. The
  reflection principle allows us to holomorphically continue
  $p(y)$ through this selected component to the double of $\cal
  M$. Therefore we can talk of the critical points of $p(y)$ on
  $(\partial {\cal M})_*$. When the argument $y$ passes through a
  simple critical point the value $p(y)$ reverses the direction
  of its movement on the circle. So there should be even number
  of critical points (counted with multiplicities) on the
  selected boundary component.

\subsection{Construction 1 (no boundary critical points)}
\label{Construct1} Using otherwise a composition with suitable
linear-fractional map, we suppose that the circle $p((\partial {\cal
M})_*)$ is the boundary of the unitary disc \be \label{Udisc}
\mbox{\sets U}:=\{p\in\mbox{\sets C}:\quad |p|\le 1\}, \ee and a
small annular vicinity of the selected boundary component is mapped
to the exterior of the unit disc. We define the mapping of a
disjoint union ${\cal M}\sqcup\mbox{\sets U}$ to a sphere:
\begin{equation}
\tilde{p}(y):=\left\{
\begin{array}{ll}
p(y),& y\in {\cal M},\\
L(y^d),& y\in \mbox{\sets U},
\end{array}
\right.
\end{equation}
where integer $d>0$ is the degree of the mapping
$p:$  $(\partial {\cal M})_*\to$ $\partial${\sets U}
and $L(y)$ is (at the moment arbitrary) linear fractional mapping
keeping the unitary disc (\ref{Udisc}). The choice of $L(\cdot)$
will be done later to simplify the arising combinatorial analysis.
\index{winding number}

  Now we fill in the hole in ${\cal M}$ by the unit disc,
  identifying the points of $(\partial {\cal M})_*$ and the
  points of $\partial${\sets U} with the same value of
  $\tilde{p}$ (there are $d$ ways to do so). The holomorphic
  mapping $\tilde{p}(y)$ of the new Riemann surface ${\cal
  M}\cup${\sets U} to the sphere will have exactly one additional
  critical point of multiplicity $d-1$ at the center of the glued
  disc.

\subsection{Construction 2 (two boundary critical points)}\label{Construct2}
  Let again $p(y)$ be a holomorphic mapping of a bounded Riemann surface $\cal
  M$ to the sphere with  selected boundary component $(\partial {\cal M})_*$
  being mapped to the boundary of the unit disc {\sets U}. Now the mapping
  $p(y)$  has two simple critical points on the selected boundary component (the
  case of coinciding critical values is not excluded). Those two points divide
  the oval $(\partial {\cal M})_*$ into two segments: $(\partial {\cal M})_*^+$
  and $(\partial {\cal M})_*^-$. We define two positive integers, {\it partial winding numbers}
  $d^\pm$ as follows. As the point $y$ moves round the selected oval in the
  positive direction, the increment of $\arg ~p(y)$ on the segment $(\partial
  {\cal M})_*^+$ is $2\pi d^+-\phi$, $0<\phi\le2\pi$, and the decrement on the
  segment $(\partial {\cal M})_*^-$ is $2\pi d^-- \phi$. We are going to modify
  the Riemann surface $\cal M$, sewing down one segment of $(\partial {\cal
  M})_*$ to the other and filling the remaining hole (if any) with the patch
  $\mbox{\sets U}$.

\begin{figure}[h]
\begin{picture}(170,50)
\thicklines
\put(30,40){\oval(40,10)}
\put(5,40){\oval(10,60)[br]}
\put(55,40){\oval(10,60)[bl]}
\thinlines
\put(30,10){\oval(50,13)[b]}
\put(20,35){\vector(-1,0){3}}
\put(40,45){\vector(1,0){3}}
\put(30,34){$*$}
\put(30,44){$*$}
\put(30,10){$\cal M$}
\put(10,48){$(\partial{\cal M})_*$}
\put(75,23){\vector(1,0){20}}
\put(80,25){$p(y)$}

\put(120,5){\psfig{figure=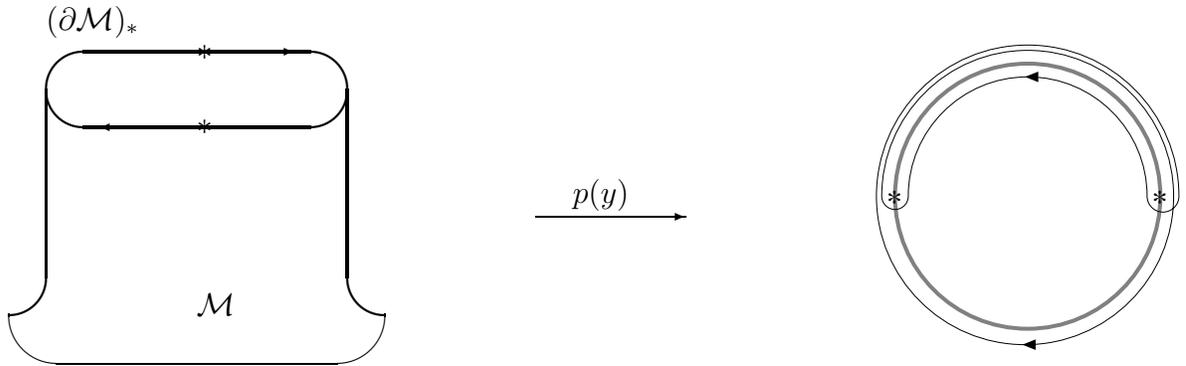}}
\end{picture}

\caption[] {\normalsize Mapping of the boundary component
$(\partial{\cal M})_*$ with two simple branch points $*$ on it and
partial winding numbers $d^+=1, d^-=2$.} \label{BdryCritPts}
\end{figure}

We define the mapping from the disjoint union ${\cal M}\sqcup\mbox{\sets U}$
to the sphere:
\begin{equation}
\tilde{p}(y):=\left\{
\begin{array}{ll}
p(y),& y\in {\cal M},\\
L(y^{d^--d^+}),& y\in \mbox{\sets U},\\
\end{array}
\right.
\end{equation}
where $L(\cdot)$ is a linear fractional mapping
keeping the unitary disc (\ref{Udisc}) invariant.

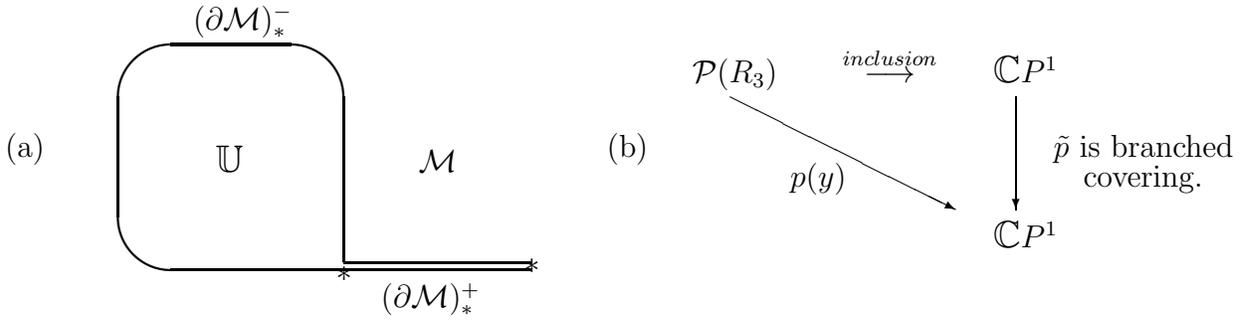
\begin{figure}
\begin{picture}(170,50)
\thicklines
\put(30,25){\oval(30,30)[t]}
\put(30,25){\oval(30,30)[bl]}
\put(30,10){\line(1,0){40}}
\put(45,11){\line(1,0){25}}
\put(45,25){\line(0,-1){14}}
\thinlines
\put(69,9.5){$*$}
\put(44,8.5){$*$}

\put(0,25){(a)}
\put(80,25){(b)}

\put(28,23){{\sets U}}
\put(55,23){$\cal M$}

\put(25,42){$(\partial{\cal M})_*^-$}
\put(50,5){$(\partial{\cal M})_*^+$}
\put(80,-5){
\begin{picture}(70,40)
\put(10,40){${\cal P}(R_3)$}
\put(30,40){$\stackrel{inclusion}\longrightarrow$}
\put(50,40){$\mbox{\sets C}P^1$}

\put(15,38){\vector(2,-1){30}}
\put(23,26){$p(y)$}
\put(53,38){\vector(0,-1){15}}
\put(58,30){$\tilde{p}$ is branched}
\put(62,26){covering.}
\put(50,18){$\mbox{\sets C}P^1$}
\end{picture}}
\end{picture}

\caption[] {\normalsize (a) Filling in the hole bounded by $(\partial{\cal M})_*$\hfill
~~~~~~~~(b) Splitting the mapping $p(y)$}
\label{DiscNCut}
\end{figure}

 Without loss of the generality we suppose that $d^-\ge d^+$. We  sew
 $(\partial {\cal M})_*^+$ to a part of $(\partial {\cal M})_*^-$, starting from
 one of the boundary critical points and consecutively identifying the points
 of the boundary oval with the same value of $p(y)$. If $d^-=d^+$ the hole disappears,
 otherwise we identify the remnant of $(\partial {\cal M})_*^-$ with the boundary of
 {\sets U} gluing points with the same value of $\tilde{p}(y)$ as shown in the Fig.
 \ref{DiscNCut}(a).

 The holomorphic mapping $\tilde{p}(y)$ of the modified Riemann surface
to the sphere will have an additional critical point of multiplicity $d^--d^+$
at the center of the artificially attached disc.
When $d^->d^+$ a simple critical point in the place of one of the old boundary critical
points arises. When $d^-=d^+$ no additional critical points arises.

\section{Proof for the Cases ${\cal A,B}1$}
\label{AB1Proof}
\subsection{Eigenfunction gives Pair of Pants}

We already know that every antisymmetric eigenfunction of integral equation PS-3
generates the mapping $p(y)$ from the pants ${\cal P}(R_3)$ to the sphere. The
boundary ovals of the pants are mapped to three circles specified in Theorem
\ref{SemiMainTh} and the function $p(y)$ may have either ({\bf a}) no critical points,
({\bf b}) one simple critical point inside the pants, ({\bf c}) two boundary simple critical
points or ({\bf d}) one double critical point on the boundary. The first two
possibilities will be considered in Sect. \ref{InnerBranch} and the other two --- in the
Sect. \ref{BoundaryBranch}

\subsubsection{No critical Points on the Boundary of Pants}\label{InnerBranch}
\paragraph{Branched Covering of a Sphere.}

  Suppose that the point $p(y)$ winds around the corresponding circle $d_r$, $d_g$ and
  $d_b$ times when the argument $y$ runs around the 'red','green' and 'blue'
  boundary component of ${\cal P}(R_3)$ respectively. We can apply the
  construction of Sect. \ref{Construct1} and glue three discs, $\mbox{\sets U}_r$,
  $\mbox{\sets U}_g$, $\mbox{\sets U}_b$ to the holes of the pants.
  Essentially, we have split our mapping $p(y)$ -- see the commutative
  diagram on the Fig. \ref{DiscNCut}(b). The holomorphic mapping
  $\tilde{p}$ has three or four ramification points, three of them are in
  the artificially glued discs and the fourth (if any) is inherited from the
  projective structure.

  Applying the Riemann--Hurwitz formula we get:

\begin{equation}
\label{RHformula}
\begin{array}{ll}
d_r+d_g+d_b=2N, & p ~ {\rm is~branched},    \\
d_r+d_g+d_b=2N+1, & p ~ {\rm is~unbranched},    \\
\end{array}
\quad N:=\deg\tilde{p}.
\end{equation}

\paragraph{Intersection of Circles.}

\begin{lmm}
In case ${\cal A}$ the required projective structure $p(y)$ with a critical point
inside the pants may exist only if the spectral parameter
$1<\lambda<2$ (i.e. when the circle $C$ does not intersect two other circles
$\varepsilon^{\pm1}\hat{\mbox{\sets R}}$).
The structure without branching does not exist for any
$\lambda$.
\end{lmm}
P~r~o~o~f:~~~
  {\bf 1.} We know that the point $0$ lies in the intersection of two
  circles: $\varepsilon\hat{\mbox{\sets R}}$ and
  $\varepsilon^2\hat{\mbox{\sets R}}$. The total number
  $\sharp\{\tilde{p}^{-1}(0)\}$ of the pre-images of this points
  (counting the multiplicities)  is $N$  and cannot be less than
  $d_b+d_g$ -- the number of pre-images on the blue and green
  boundary oval of the pants. Comparing this to
  (\ref{RHformula}) we get $d_r\ge N$ which is only possible when

\begin{equation}
\label{k1k2kRelation}
d_r=d_g+d_b=N.
\end{equation}

  Assuming that the circle $C$ intersects any of the circles
  $\varepsilon^{\pm1}\hat{\mbox{\sets R}}$ we repeat the above
  argument for the intersection point and arrive at the
  conclusion $d_b=d_r+d_g=N$ or $d_g=d_r+d_b=N$ incompatible with
  already established (\ref{k1k2kRelation}).

  {\bf 2.} For the unbranched structure $p(y)$ the established inequalities
$d_b+d_g\le N$ and $d_r\le N$ contradict the  Riemann-Hurwitz
formula (\ref{RHformula}). ~~~\bl

The above arguments may be applied to the case ${\cal B}1$ as well.
Taking into account that the circles $C$ and $\hat{\mbox{\sets R}}$ always intersect
we arrive at
\begin{lmm}
In case ${\cal B}1$ the mapping $p(y)$ (if any) will have a
boundary critical point.
\end{lmm}

\begin{figure}[ht!]
\begin{picture}(170,55)
\put(0,0){\psfig{figure=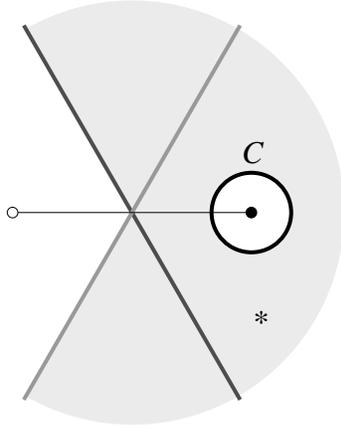}}
\end{picture}
\caption[] {\normalsize Shaded area is the image of pants in case $\cal A$, index $d_b>0$.}
\label{BelyiMap}
\end{figure}

\paragraph{Image of Pants.}
  Let us investigate where the artificially glued discs in case $\cal A$ are
  mapped to. Suppose for instance that the disc $\mbox{\sets
  U}_r$ is mapped to the exterior of the circle $C$. The point
  $0$ will be covered then at least $d_r+d_g+d_b=2N$ times which
  is impossible. The discs $\mbox{\sets U}_g$ and $\mbox{\sets
  U}_b$  are mapped to the left of the lines
  $\varepsilon{\mbox{\sets R}}$ and $\varepsilon^2{\mbox{\sets
  R}}$ respectively, otherwise points from the interior of the
  circle $C$ will be covered more that $N$ times.  The image of
  the pair of pants ${\cal P}(R_3)$ is shown on the
  Fig. \ref{BelyiMap}.

  We use the ambiguity in the construction of gluing the discs to the pants and
  require  that the critical values of $\tilde{p}$ in the discs $\mbox{\sets
  U}_g$, $\mbox{\sets U}_b$ coincide. Now the branched covering $\tilde{p}$ has
  only three different branch points shown as $\bullet$, $\circ$ and $*$ on the
  Fig. \ref{BelyiMap}. The branching type at  $\bullet$ is the cycle of
  length $N$; at the point $\circ$ there are cycles of lengths $d_g$ and $d_b$;
  and the branch point $*$ is simple. The coverings with three branch points are
  called Belyi maps and are described by certain graphs known as Grothendieck's
  {\it "Dessins d'Enfants"}. In our case the {\it dessin} is the lifting of the
  segment connecting white and black branch points: $\Gamma:=\tilde{p}^{-
  1}[\bullet,\circ]$.

\begin{figure}[h!]
\psfig{figure=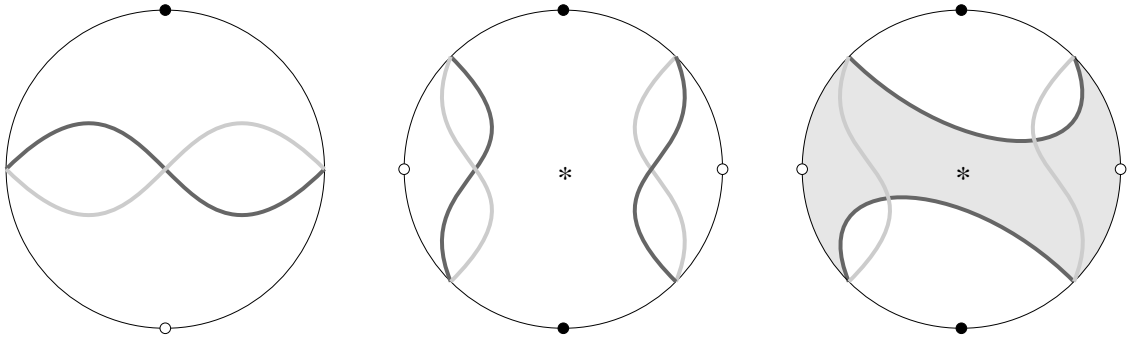}
\caption[] {\normalsize (a) Simple cell ($N-2$ copies)~~ (b) Double cell (1 copy)~~ (c) Impossible double cell}
\label{TwoCells}
\end{figure}

\paragraph{Combinatorial Analysis of Dessins.}
  There is exactly one critical point of $\tilde{p}$ over the branch point $*$.
  Hence, the compliment to the graph $\Gamma$
  on the upper sphere of the diagram on the Fig.
  \ref{DiscNCut}(b) contains exactly one cell mapped $2-1$ to the lower sphere.
  All the rest components of the compliment are mapped $1-1$.
  Two types of cells are shown
  in the Fig.  \ref{TwoCells} (a) and (b), the
  lifting of the red circle is not shown to simplify the
  pictures. The branch point $*=:h_1+ih_2$ should lie in the intersection of
  two annuli $\alpha$ and $\overline{\alpha}$ otherwise
  the discs $\mbox{\sets U}_g$, $\mbox{\sets U}_b$ glued to different boundary components of
  our pants will intersect: the hypothetical case when the
  branch point of $p(y)$ belongs to one annulus but does not
  belong to the other is shown in the Fig. \ref{TwoCells} (c).

\begin{figure}[ht!]
\psfig{figure=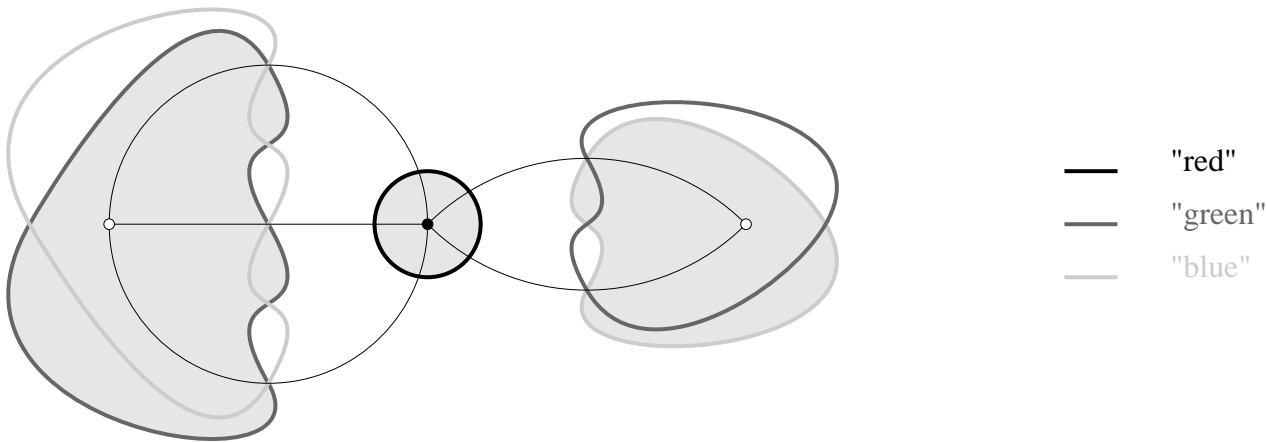} \caption[] {\normalsize Dessin for
$d_g=3,d_b=2$; the pre-image of the branch point $*$ is at infinity}
\label{Assembly1}
\end{figure}

  The cells from the Fig. \ref{TwoCells} (a), (b) may be assembled in
  a unique way shown in the Fig. \ref{Assembly1}. The pants are
  colored in  white, three artificially sewn discs are
  shaded. Essentially this picture shows us how to sew together
  the patches bounded by our three circles $C$,
  $\varepsilon^{\pm1}\hat{\mbox{\sets R}}$ to get the pants
  conformally equivalent to ${\cal P}(R_3)$. As a result of the surgery procedure we get
  the pants ${\cal PA}_1(\lambda, h_1, h_2|d_g,d_b)$. Changing the superscript
  of the projective structure $p^\pm(y)$ gives us the change of sign for the eigenfunction
  $u(x)$ and the reflection of the pants ${\cal PA}_1(\dots)$ in the unit circle
  $\partial${\sets U}. This is why we consider only the pants with $|h_1+ih_2|\ge1$.

\subsubsection{Boundary critical Points}
\label{BoundaryBranch}
First of all we consider the stable case of two simple critical points on the boundary
oval. At the moment we do not know the color of this oval and we use the
'nicknames' $\{1,2,3\}$  for the set of colours $\{r,g,b\}$ so that the critical
points will be on the oval 3.

\paragraph{Branched Covering of a Sphere.}
 The usage of both constructions from Sect.  \ref{Construct} allows us to
 include the pants ${\cal P}(R_3)$ to the sphere attaching two
 discs $\mbox{\sets U}_1$ and $\mbox{\sets U}_2$ to the first two ovals and
 collapsing the boundary of the third oval and sewing in the third disc
 $\mbox{\sets U}_3$ if necessary. Positive integers arising in those
 constructions we denote as $d_1$, $d_2$, $d_3=|d_3^--d_3^+|$ respectively.

  We arrive at the branched covering $\tilde{p}$ of the diagram
  on the Fig \ref{DiscNCut}(b). This mapping has either two or four
  critical points. Two of them are in the centers of the discs
  $\mbox{\sets U}_1$ and $\mbox{\sets U}_2$, another two
  arise only when $d_3>0$ -- the center of $\mbox{\sets U}_3$
  and one of the boundary critical points of the mapping $p(y)$.
  The multiplicities of those critical points are respectively
  $d_1-1$, $d_2-1$, $d_3-1$, 1. Riemann-Hurwitz formula for this covering reads
 \be
 \label{RH2}
  d_1+d_2+d_3=2N,
  \qquad N:=\deg\tilde{p}.
 \ee

\begin{lmm}
The images of the ovals 1 and 2 do not intersect.
\label{ovals12}
\end{lmm}

P~r~o~o~f. Suppose the inverse is true and a point $Pt$ lies in the
intersection of images of the first two ovals. Then
$N\ge\sharp\tilde{p}^{-1}(Pt)\ge d_1+d_2$. Any of the critical
points from the third oval has at least $d_3+1\le N$ pre-images counting
multiplicities. The last two inequalities contradict (\ref{RH2}).
~~\bl

{\bf Corollaries.}

{\bf 1.} In case $\cal A$ the critical points of $p(y)$
lie either on the blue or on the green boundary of pants.
(Two circles $\varepsilon^{\pm1}\hat{\mbox{\sets R}}$ intersect)

{\bf 2.} In case ${\cal B}1$ the critical points of $p(y)$ lie on the
blue boundary of pants. (Two circles $C$ and $\hat{\mbox{\sets R}}$
intersect)

{\bf 3.} In both cases the required function may only exist when
$\mu\in(\frac12,1)$, or equivalently $\lambda\in(1,2)$.
(Otherwise the circles  $C$ and $\varepsilon^{\pm1}\hat{\mbox{\sets R}}$ intersect)

To save space, further proof will be given for the case $\cal A$ only when both
critical points lie on the blue oval. The omitted cases require no extra technique.
Now the notations $\mbox{\sets U}_r$, $\mbox{\sets U}_g$, $\mbox{\sets U}_b$,
$d_r$, $d_g$, $d_b$ have the obvious meaning.

\paragraph{Image of Pants.}
\begin{lmm}
The image $p({\cal P}(R_3))$ of the pants is the union $\alpha\cup\bar\alpha$
when $d_b>0$ or the annulus $\alpha$ when $d_b=0$.
\end{lmm}
P~r~o~o~f:
Essentially, we have established the equalities
$$d_g+d_b=d_r=N$$ setting $Pt=0$ in the proof of Lemma \ref{ovals12}.
Repeating the arguments of the same title paragraph of the
Sect. \ref{InnerBranch} we conclude that the disc $\tilde{p}(\mbox{\sets U}_r)$ fills
the interior of $C$, the disc $\mbox{\sets U}_g$ is mapped to the left of
the line $\varepsilon{\mbox{\sets R}}$ and the disc $\mbox{\sets U}_b$ (if any) is
mapped to the left of $\varepsilon^2{\mbox{\sets R}}$.
So the sector $\{\frac{2\pi}3\le\arg p\le\frac{4\pi}3 \}$
is covered $d_g+d_b=N$ times by the artificially inserted discs.
The disc $\mbox{\sets U}_g$ covers the half-plane to the left of
the line $\varepsilon{\mbox{\sets R}}$
exactly $d_g$ times, the latter number is $N$ when $d_b=0$.

{\bf Corollary} ~~ Both critical values of $p(y)$ lie on the ray
$-\varepsilon^2(0,\infty)$.

\begin{figure}[h!]
\psfig{figure=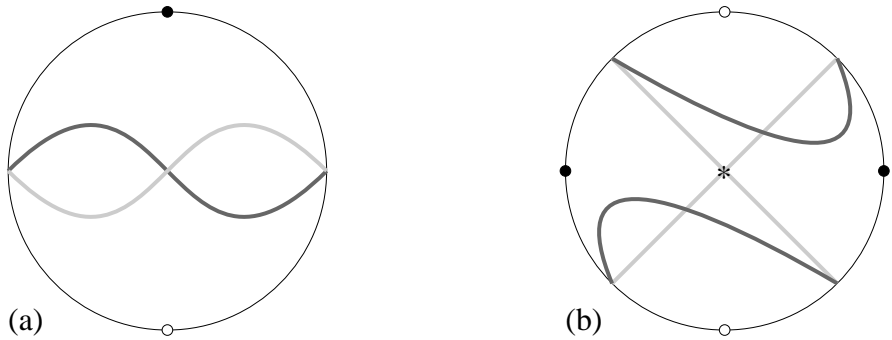}
\caption[] {\normalsize (a) Simple cell ($N-2$ copies)~~ (b) Double cell (1 copy)}
\label{GraphCells2}
\end{figure}

\paragraph{Dessin d'Enfants.}
  Again, we put the critical values of $\tilde{p}$ in the discs $\mbox{\sets
  U}_g$, $\mbox{\sets U}_b$ to the same point $\circ$ (see Fig.
  \ref{BelyiMap}). The only difference from the Sect. \ref{InnerBranch}:
  now the inherited from the pants branch point $*$ (if $d_b>0$) lies on the ray $-
  \varepsilon^2(0,\infty)$. We introduce the Grothendieck's {\it Dessin} as the
  lifting of the segment connecting white and black branch points:
  $\Gamma:=\tilde{p}^{- 1}[\bullet,\circ]$. The compliment to $\Gamma$
  is composed of cells shown in the Fig. \ref{GraphCells2}. In the assembly the
  double cell may be used only once and only when $d_b>0$.

\begin{figure}[ht!]
\psfig{figure=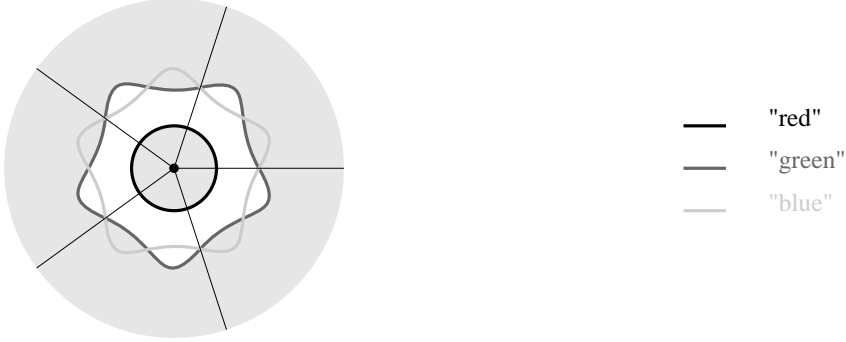} 
\caption[] {\normalsize Degenerate
Dessin for $d_g=5,d_b=0$; the pre-image of the branch point $\circ$
is at infinity} \label{Assembly2}
\end{figure}

  Given the winding numbers $d_g,d_b$, the cells from the Fig. \ref{GraphCells2} 
  (a), (b) may be attached to each other in a unique way. When $d_b>0$ the 
  {\it Dessin} $\Gamma$ has the same combinatorial structure as in Fig. 
  \ref{Assembly1}. Of course one  has to replace the old cells by those shown in 
  Fig. \ref{GraphCells2}. Shown in the Fig. \ref{Assembly2} is the assembly of 
  cells for $d_b=0$, $d_g=5$. The pants are colored in  white, two artificially 
  inserted discs are shaded.  As a result of the surgery procedure we get the 
  pants ${\cal PA}_2(\lambda, h_1, h_2|d_g,d_b)$ with positive reals $h_1,h_2$ 
  determined by the critical values of $p(y)$. To discern the pair of pants 
  ${\cal PA}_2(\dots)$ from its reflection in the unit circle we consider the 
  restriction $h_1h_2\ge1$.

\paragraph{Junction of critical points.}
To study the remaining case when the boundary critical points of projective 
strucure merge, one has to apply the limit case of the Construction 2. In this 
way one arrives at the unstable membranes ${\cal PA}_{12}$ and ${\cal PA}_{13}$. 
To save space we omit the details.

\subsection{Pair of Pants corresponds to Eigenfunction}
Let the pair of pants (\ref{R3Pants}) associated to the functional
parameter of the type ${\cal Q=A,B}1$ integral equation is
conformally equivalent to the pants ${\cal PQ}(\lambda, \dots)$.
This exactly means that there exists a conformal mapping $p(y)$
from ${\cal P}(R_3)$ to ${\cal PQ}(\lambda, \dots)$
respecting the colors of the boundary ovals. This mapping
is unique since the conformal self-mapping of pants keeping all boundary ovals invariant 
is trivial.  The mapping $p(y)$ has one simple critical point inside the pants 
(for the membrane ${\cal PA}_1(\dots)$) or two simple boundary points (for  ${\cal PA}_2(\dots)$,
${\cal PA}_3(\dots)$, ${\cal PB}1(\dots)$) or a double boundary
critical point (for ${\cal PA}_{12}(\dots)$, ${\cal PA}_{13}(\dots)$). 
Moreover, $p(y)$ maps the boundary components of ${\cal P}(R_3)$ to the 
circles specified by Theorem \ref{SemiMainTh}. Hence, given $p(y)$ one can consecutively
restore: two projective structures $p^\pm(y)$, the solution $W(y)$
of Riemann monodromy problem and the eigenfunction $u(x)$. Combining
the formulae (\ref{Antisymm}), (\ref{UniU}) we obtain the top of the
reconstruction formulae in (\ref{UviaP}).

\section{Proof for the Case ${\cal B}2$}
\label{B2Proof}
\subsection{Eigenfunction gives Pair of Pants}
\begin{figure}[ht!]
\psfig{figure=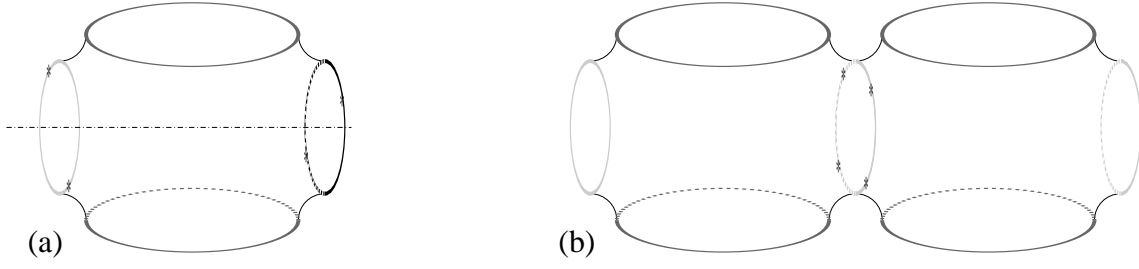} 
\caption[] {\normalsize (a) The surface ${\cal P}_2$ is the double cover of pants ${\cal P}$;\\ 
(b) ${\cal P}_4$ is the reflection of ${\cal P}_2$ in the blue boundary oval }
\label{AuxSurface}
\end{figure}

Any antisymmetric eigenfunction of the integral equation PS-3 
generates the mapping $p(y)$ from the pants ${\cal P}(R_3)$ to the sphere.
The principal difference of this case from the one considered in Sect. 
\ref{AB1Proof} lies in the two-valuedness of the function $p(y)$ in the pants. 
To reflect this phenomenon we consider the two sheeted 
unbranched cover ${\cal P}_2\to$ ${\cal P}(R_3)$ with trivial monodromy around 
the green boundary oval. This new surface is a sphere with four holes, each 
boundary inherits the color of the oval it covers -- see Fig. \ref{AuxSurface}(a). The mapping $p(y)$ is lifted 
to the single-valued mapping $p_2:~~$ ${\cal P}_2\to$ $\mbox{\sets C}P^1$ 
satisfying the equivariance condition:

\be
\label{EqVar}
p_2\Xi=\chi({\bf DD}_1)p_2,
\ee
where $\Xi$ is the cover transformation (change of sheets) of
${\cal P}_2$.

\vspace{1mm}

Now we can complete the p~r~o~o~f~ o~f~ t~h~e~ L~e~m~m~a
\ref{SpLocus}.\\ Suppose that there exists the required function
$p_2(y)$ in the case ${\cal B}21$ and $\mu>1$. We show that all
possible locations of the critical points of this function lead to
the contradiction: ({\bf a}) $p_2$ has no boundary critical points;
({\bf b}) all critical points lie on the green ovals of ${\cal P}_2$;
({\bf c}) $p_2$ has at least one critical point on a blue oval.

({\bf a}) We use the Construction 1 and attach four discs to the surface
${\cal P}_2$. For the arising ramified covering $\tilde{p}_2$ the Riemann-Hurwitz
formula reads
$$
\begin{array}{ll}
2d_g+d_b+d_b'=2N, & p ~ {\rm is~branched},  \\
2d_g+d_b+d_b'=2N+2, & p ~ {\rm is~unbranched},  \\
\end{array}
\quad N:=\deg\tilde{p}_2,
$$
where $d_g$ is the winding number of $p_2(y)$ for each of green ovals; $d_b$ and
$d_b'$ are the winding numbers for two blue ovals of the surface ${\cal P}_2$. The point
$0\in$ $\hat{\mbox{\sets R}}\cap$ $\hat{\varepsilon\mbox{\sets R}}$ is covered
at least $d_g+d_b+d_b'\le N$ times. This agrees the previous formula only if
$d_g\ge N$. But now $d_b=d_b'=0$ which is impossible.

({\bf b}) Now each of two green ovals has two boundary critical points of $p_2(y)$. We
use both Constructions and eliminate all holes in ${\cal P}_2$ attaching two
discs to the blue ovals and possibly two more  discs to the green ovals of the
surface. The Riemann-Hurwitz formula for the arising ramified covering
$\tilde{p}_2$ reads
$$
2d_g+d_b+d_b'=2N,
\quad N:=\deg\tilde{p}_2,
$$
where $d_g\ge0$ is the difference of the partial winding numbers for each of the 
green ovals of ${\cal P}_2$. Further argument is exactly as in the previous 
paragraph.

({\bf c})
We claim that in this case there are exactly four critical points of $p_2$ on a 
blue oval of the surface ${\cal P}_2$. Indeed, given a critical point $Pt$, $\Xi 
Pt$ is also a critical point because of the equivariance (\ref{EqVar}). When 
$\mu>1$, the mapping $\chi({\bf DD}_1)$ conserves the orientation of the real axis. 
This means that those two critical points are of the same type (say, local 
maxima of the real value $p_2$). Hence, $Pt$ and $\Xi Pt$ are separated by the 
critical points of the opposite type (local minima in our case). There cannot be 
more that four critical points of the function $p_2$ on the double cover of 
pants ${\cal P}(R_3)$, so we have listed them all.

Let us consider the {\it double} of the surface ${\cal P}_2$
and cut it along all boundary ovals of  ${\cal P}_2$, but
the blue
oval containing all critical point of $p_2(y)$. This new surface, ${\cal P}_4$, is a
sphere with six holes shown in the Fig \ref{AuxSurface}(b), four boundary ovals are 
green and two are blue.  The
reflection principle allows to continue analytically $p_2(y)$ to the mapping
$p_4(y)$ of the entire surface ${\cal P}_4$ to the sphere. This continuation has
four inner critical points and no boundary critical points. It maps both blue
ovals to $\hat{\mbox{\sets R}}$ and four green ovals to the circles
 $\varepsilon^{\pm1}\hat{\mbox{\sets R}}$,  
 $\chi({\bf DD_1)})\varepsilon^{\pm1}\hat{\mbox{\sets R}}$. The
usage of Construction 1 allows to fill in all the holes of ${\cal P}_4$. The
Riemann-Hurwitz formula for the arising ramified covering $\tilde{p}_4$
reads
$$
2d_g+d_b=N,
\quad N:=\deg\tilde{p}_4,
$$
where the numbers $d_g$, $d_b$ have the obvious meaning. The point $0$ is covered at least
$2d_g+2d_b>N$ times which is impossible.
~~~\bl

\paragraph{The location of the critical points} 
of the mapping $p_2(y)$ is given by the following lemma.

\begin{lmm}
The mapping $p_2(y)$ has exactly two boundary critical points  
on each of the non-green ovals of the surface ${\cal P}_2$.
\label{p2Branch}
\end{lmm}

P~r~o~o~f. The mapping $\chi({\bf DD}_1)$ changes the orientation of the circles
$C$ and $\hat{\mbox{\sets R}}$, when $\mu\in(0,1)$. When the point $y$ runs along
the blue or red oval of ${\cal P}_2$, the value $p_2(y)$ changes the orientation
of its motion at least twice due to (\ref{EqVar}). This means that the argument
$y$ comes through at least two boundary critical points. Since the mapping
$p(y)$ has at most two boundary critical points in pants, the lifted mapping $p_2(y)$ has
at most four in the double cover of pants. ~~~\bl

\paragraph{Branched covering of the sphere.}
The mapping $p_2(y)$ from ${\cal P}_2$ to the sphere has equal
winding numbers $d=d_g$ defined in Sect. \ref{Construct1} on both green ovals.
On each of non-green ovals, $p_2(y)$ has zero index $d:=d^--d^+=0$ introduced in 
Sect. /ref{Construct1}.
Both statements are simple consequences of the equivariance condition (\ref{EqVar}).
Applying Constructions 1 and 2 to the mapping $p_2$ defined on ${\cal P}_2$,
we get a ramified covering $\tilde{p}_2$ with two critical points, both of multiplicity
$d_g-1$.

\paragraph{Image of the Surface.}
The Riemann-Hurwitz formula for the ramified covering $\tilde{p}_2$ reads
$d_g=N:=\deg \tilde{p}_2$. It is easily seen that two discs attached to
the green ovals of ${\cal P}_2$ are mapped to the left of the line $\varepsilon$~{\sets R}
and to the interior of the circle $\chi({\bf DD}_1)$ $\varepsilon\hat{\mbox{\sets R}}$.
Therefore, the surface ${\cal P}_2$ is conformally equivalent to the
closure of the annulus $d_g\cdot\beta$ with two slots in it.

The involution $\Xi$ of ${\cal P}_2$ (the interchange of sheets) induces the involution 
of the multisheeted annulus. The latter involution is the lifting of $\chi({\bf 
DD}_1)$ to $d_g\cdot\beta$ and is given by the formula (\ref{mBetaInv}). The 
slots of $d_g\cdot\beta$ are invariant with respect to $\Xi$ and therefore pass 
through the fixed points $\mu^{-1}+r$ and $\mu^{-1}+r\exp{i\pi m}$ of the 
involution. The red slots are projected to the circle $C$, the blue slots are 
projected to the real line. Given in Tab. \ref{SlotsDef} inequalities for the 
parameters $h_1,h_2$ specifying the endpoints of the slots allow us to relate any 
given antisymmetric eigenfunction to exactly one picture.

A by-product of the explicit description of the image of the pants is the following  
\begin{lmm}
\label{PpmDifferent}
In antisymmetric case ${\cal B}2$ two structures $p^\pm(y)$ are different.
\end{lmm}

P~r~o~o~f.~ Suppose the opposite, that is 
\be
p(y)\overline{p(\bar{y})}\equiv 1
\label{StructEqual}
\ee
for $p(y)$ satisfying the conditions of theorem \ref{SemiMainTh}.

In case ${\cal B}21$ the value $p(a)\in${\sets R} when  $a\in \{a_1,a_2\}$
is the endpoint of the blue slot. From (\ref{StructEqual})
it immediately follows that $p(a)=\pm1$. But the image of pants
$p({\cal P})=\beta$ avoids both points $\pm1$. 

In case ${\cal B}22$ the value $p(a\pm i0)\in C=\{p=\chi({\bf DD_1})\bar{p}\}$  
when  $a\in \{a_1,a_2\}$ is the endpoint of the red slot. From (\ref{StructEqual})
and the jump relationship (\ref{Jump}) on $[a_1,a_2]$
it follows that $p(a\pm i0)=\pm1$. Again, the image of pants
$p({\cal P})$ avoids both points $\pm1$. 

In case ${\cal B}23$ any of the above two arguments is applicable. ~~~\bl
 
{\bf Corollary~ }
{\it In case ${\cal B}2$ any eigenvalue corresponds to no more than one 
antisymmetric eigenfunction.}

P~r~o~o~f. Suppose that three meromorphic functions $p_s(y)$, $s=1,2,3$, in the pants satisfy 
the conditions of theorem \ref{SemiMainTh} and no two of them are identical.
From the second part of the proof of the theorem \ref{J0neq0} we know
that all three values $p_s(y)$ are different at any point $y$.
We consider the following differential form on the Riemann surface $M$:
$$
\omega:=dp_1
\left(
\frac1{p_1-p_2}-\frac1{p_1-p_3}
\right).
$$

This form $\omega$ is the infinitesimal form of cross ratio and it is invariant 
under the same linear-fractional transformations of all three functions $p_s$. 
Therefore $\omega$ is well defined on the entire Riemann surface $M$. Using 
local coordinates on $M$, it's easy to check that the form is holomorphic and 
$(\omega)={\sf D}(p_1)$. Any holomorphic differential on the genus 2 surface has 
two zeroes which are interchanged by the hyperelliptic involution of $M$. 
According to Lemma \ref{p2Branch}, the branching divisor of $p_1(y)$ is 
different as it has a branchpoint on each of  non-green ovals of the pants $\cal 
P$. 

Therefore, two of our functions $p_s(y)$ coincide. 
Moreover,   from Lemma \ref{PpmDifferent} it follows that either all three functions 
are identical, or one of them is the antisymmetrization of the other: 
$p_2(y)=1/\overline{p_1(\bar{y})}$  and $p_3=p_1$. ~~~\bl

\subsection{Pair of Pants corresponds to Eigenfunction}
Let the pair of pants ${\cal PB}2s(\lambda,h_1,h_2|m)$ be conformally equivalent to the 
pair of pants ${\cal P}(R_3)$ associated to type ${\cal B}2s$, $s=1,2,3$, integral equation.
This exactly means that there exists respecting the colors of the boundary ovals
equivariant conformal mapping 
$p_2(y)$ from the double cover ${\cal P}_2(R_3)$ of pants to the closure of the 
multisheeted annulus $m\cdot\beta$ with two slots $E_s^1(h_1)$ and $E_s^2(h_2)$ in it. 
We represent the double cover ${\cal P}_2(R_3)$ as two copies of pants ${\cal P}(R_3)$
cut along the segment $[-1,1]\cap[a_1,a_2]$ and attached one to the other.
The restriction of $p_2(y)$ to one of such copies gives the function $p(y)$ satisfying 
all the assumptions of the theorem \ref{SemiMainTh}. The antisymmetric eigenfunction of 
the integral equation now may be reconstructed via the known procedure
which gives the formulae (\ref{UviaP}). Moreover, from Lemma \ref{PpmDifferent}
we have learned that the invariant $J_0\neq0$ for antisymmetric solutions in case
${\cal B}2$. So we can use the alternative formula (\ref{UfromP})
to reconstruct the eigenfunction $u(x)$ when $x\in[-1,1]\setminus[a_1,a_2]$.
For $x\in[-1,1]\cap[a_1,a_2]$ we can also use the following formula:
$$
\frac{{\rm Im}~p(y^+)}{|p(y^+)-\mu|^2+1-\mu^2},
\quad y^+:=R_3(x+i0),
\quad x\in[-1,1]\cap[a_1,a_2].
$$

The only nuisance here consists in possible non-uniqueness of the mapping $p(y)$.
Indeed, when two of the boundary ovals of pants
have the same color (blue in case ${\cal B}21$ of red in case ${\cal B}22$),
the pants may admit conformal involution interchanging the ovals of the same color.
Such pants fill in a codimension one manifold in the corresponding moduli space.  
The Corollary to Lemma \ref{PpmDifferent} nevertheless guarantees the uniqueness of 
the antisymmetric eigenfunction for the given membrane ${\cal PB}2s$: the composition of 
$p(y)$ with the conformal automorphism of pants coincides with either $p(y)$ 
or its antsymmetrization $1/\overline{p(\bar{y})}$.

\section{Conclusion}
Similar analysis based on the geometry and combinatorics may be applied to obtain the
representations of the solutions of PS-3 integral equation in all the dropped cases.
Much of the techniques used is helpful for the study of other integral equations with
rational low degree kernels.

\vspace{5mm}
\parbox{8cm}
{\it
119991 Russia, Moscow GSP-1,\\
ul. Gubkina 8,\\
Institute for Numerical Mathematics,\\
Russian Academy of Sciences}


\begin{thebibliography}{}




  \bibitem{Bog0} Bogatyrev A.B. {\it ~The discrete spectrum of the
  problem with a pair of Poincare-Steklov operators.}// Doklady
  RAS 358:3, (1998). See also
  Bogatyrev A.B. {\it ~The spectral properties of
  Poincare-Steklov operators.} PhD Diss., INM RAS, Moscow 1996.

\bibitem{Bog3} Bogatyrev A.B. {\it A geometric method for solving a series of
  integral PS equations} // Math. Notes, 63:3 (1998), pp. 302-310.

\bibitem{Poi} Poincare H. {\it Analyse des travaux scientifiques de Henri Poincare}
 //Acta Math 38 (1921), pp. 3-135.

  \bibitem{G} Gunning R.C. {\it Special coordinate coverings of
  Riemann surfaces}// Math. Annalen, 170(1967), pp. 67-86.

  \bibitem{Tu} Tyurin A.N. {\it On the periods of quadratic
  differentials}. //Russian Math. Surveys, 33:6 (1978).

  \bibitem{GKM} Gallo D., Kapovich M., Marden A. {\it The
  monodromy groups of Schwarzian equations on closed Riemann
  surfaces} //Ann. of Math. (2) {\bf 151}:2, 625-704 (2000).
  See also arXiv, math.CV/9511213.

  \bibitem{Hej} D.A.Hejhal {\it Monodromy groups and linerly
  polymorphic functions} //Acta Math.,135(1975), pp. 1-55.

  \bibitem{Man} R.Mandelbaum {\it Branched structures and affine
  and projective bundles on Riemann surfaces} //Trans. AMS,
  183(1973), pp. 37-58.

  \bibitem{Bog2} Bogatyrev A.B. {\it Poincare-Steklov integral
  equations and the Riemann monodromy problem}// Funct. Anal.
  Appl. 34:2 (2000), pp. 9-22.

  \bibitem{Bog1} Bogatyrev A.B. {\it PS-3 integral equations and
  projective structures on Riemann surfaces}// Sbornik:
  Mathematics 192:4 (2001), pp. 479-514.



\end{thebibliography}
\end{document}